\documentclass[aop,MSNbibl,seceqn,number,citesort,dvips]{arximspdf}

\doi{10.1214/11-AOP703}
\volume{41}
\issue{3A}
\pubyear{2013}
\firstpage{1461}
\lastpage{1512}

\makeatletter

\newtheorem{theorem}{Theorem}[section]
\newtheorem{lemma}[theorem]{Lemma}
\newtheorem{proposition}[theorem]{Proposition}
\newtheorem{corollary}[theorem]{Corollary}
\renewcommand{\c}{\mathfrak c}
\def\J{\mathbb J}
\def\Z{\mathbb Z}
\def\R{\mathbb R}
\def\I{\mathbb I}
\newcommand{\eqref}[1]{(\ref{#1})}
\newcommand{\fracd}[2]{({#1}/{#2})}
\newcommand{\fracc}[2]{{#1}/{(#2)}}
\newcommand{\fraca}[2]{{#1}/{#2}}
\makeatother

\begin{document}
\begin{frontmatter}

\title{Large deviations for the current and tagged particle in 1D
nearest-neighbor symmetric simple~exclusion}
\runtitle{LDP for current and tagged particle in 1D SSEP}

\begin{aug}
\author[A]{\fnms{Sunder} \snm{Sethuraman}\corref{}\thanksref{au1}\ead[label=e1]{sethuram@math.arizona.edu}}
\and
\author[B]{\fnms{S. R. S.} \snm{Varadhan}\thanksref{au2}\ead[label=e2]{varadhan@cims.nyu.edu}}
\runauthor{S. Sethuraman and S. R. S. Varadhan}
\affiliation{University of Arizona and New York University}
\address[A]{Department of Mathematics\\University of Arizona\\617 N.
Santa Rita Ave.\\Tucson, Arizona 85721\\USA\\\printead{e1}} 
\address[B]{Courant Institute of Mathematical Sciences\\New York
University\\212 Mercer St.\\New York, New York 10012\\USA\\\printead{e2}}
\end{aug}
\thankstext{au1}{Supported in part by NSF-0906713 and NSA Q H982301010180.}
\thankstext{au2}{Supported in part by NSF Grants DMS-09-04701 and
OISE-0730136.}

\received{\smonth{1} \syear{2011}}
\revised{\smonth{7} \syear{2011}}

%
\begin{abstract}
Laws of large numbers, starting from certain nonequilibrium measures,
have been shown for the integrated current across a bond, and a tagged
particle in one-dimensional
symmetric nearest-neighbor simple exclusion [\textit{Ann. Inst. Henri
Poincar\'e Probab. Stat.} \textbf{42} (2006) 567--577].
In this article, we
prove corresponding large deviation principles and evaluate the rate
functions, showing different growth behaviors near and far from their
zeroes which connect with results in [\textit{J. Stat. Phys.} \textbf
{136} (2009) 1--15].
\end{abstract}

%

\end{frontmatter}

\section{Introduction and results}
\label{sec0}

The one-dimensional
nearest-neighbor symmetric simple exclusion process
follows a collection of nearest-neighbor random walks on
the lattice $\Z$, each of which is equally likely to move left or right,
except in that
jumps to already occupied sites are suppressed.
More precisely, the model is a Markov process $\eta_t=
\{\eta_t(x)\dvtx  x\in\Z\}$, evolving on the configuration space $\Sigma
= \{0,1\}^\Z$, with generator
\[
(L\phi)(\eta) = (1/2)\sum_x \bigl[ \eta(x)\bigl(1-\eta(x+1)\bigr
)+\eta(x+1)\bigl(1-\eta(x)\bigr) \bigr] \bigl( \phi(\eta
^{x,x+1}) -
\phi(\eta) \bigr),
\]
where $\eta^{x,y}$, for $x\neq y$, is the configuration obtained from
$\eta$ by
exchanging the values at $x$ and $y$,
\[
\eta^{x,y}(z) =
\cases{
\eta(z), &\quad when $z\neq x,y$,\cr
\eta(x), &\quad when $ z=y$,\cr
\eta(y),&\quad when $z=x$.}
\]
A detailed
treatment can be found in Liggett~\cite{Liggett}.

As the process is ``mass conservative,'' that is, no birth or death,
one expects a family of invariant measures corresponding to particle
density. In fact, for each $\rho\in[0,1]$, the product over $\Z$
of Bernoulli measures $\nu_\rho$ which independently puts a particle
at locations $x\in\Z$ with probability $\rho$, that is,
$\nu_\rho(\eta_x=1)=1-\nu_\rho(\eta_x=0)=\rho$, are invariant.
We will denote $E_\rho$ as expectation under $\nu_\rho$.

Consider now the integrated current across the bond $(-1,0)$, and a
distinguished, or tagged particle, say initially at
the origin. Let $J_{-1,0}(t)$ and $X_t$ be the current and position of
the tagged particle at time $t$, respectively. The problem of
characterizing the asymptotic behavior of the current and tagged
particle in interacting
systems has a long history (cf. Spohn~\cite{Spohn}, Chapters 8.I, 6.II),
and was mentioned in Spitzer's seminal paper~\cite{Spitzer}.

The goal of this paper is to investigate the large deviations of
$J_{-1,0}(t)$ and $X_t$ when
the initial distribution of particles is part of a large class of
nonequilibrium measures.
Our initial motivation was to understand better laws of large numbers
(LLN) and central limit theorems (CLT) in Jara and Landim~\cite{JL}
for the current and tagged particle when the process starts from a
class of ``local equilibrium'' initial measures. It turns out that
recent formal expansions of the large deviation ``pressure'' for the
current in Derrida and Gerschenfeld~\cite{Derrida,Derrida2}
might also be recovered in such a study.

The article~\cite{JL} is a
nonequilibrium generalization of
CLTs in Arratia
\cite{Arratia}, Rost and Vares~\cite{Rost-Vares}, and De Masi and
Ferrari~\cite{DeMF}, which established ``subdiffusive'' behaviors in
the 1D
nearest-neighbor symmetric simple exclusion model. Namely, starting
under an equilibrium
$\nu_\rho$, $t^{-1/4}J_{-1,0}(t) \Rightarrow\mathrm{N}(0,\sigma
_J^2)$ and $t^{-1/4}X_t \Rightarrow\mathrm{N}(0,\sigma_X^2)$, where
$\sigma_J^2 = \sqrt{2/\pi}(1-\rho)\rho$ and $\sigma_X^2 = \sqrt
{2/\pi}(1-\rho)/\rho$.
Physically, the ``subdiffusive'' scale in the CLT is explained as being due
to ``trapping'' induced from the nearest-neighbor dynamics which
enforces a rigid ordering of particles. Recently, the CLTs were
extended to an invariance principle with respect to a fractional
Brownian motion, $\lambda^{-1/4}J_{-1,0}(\lambda t) \Rightarrow\sigma
_J fBM_{1/4}(t)$ and $\lambda^{-1/4}X_{\lambda t} \Rightarrow\sigma
_X fBM_{1/4}(t)$, in Peligrad and Sethuraman~\cite{Peligrad-S}.

We now specify the class of initial measures considered, that is,
``deterministic initial configurations'' and ``local equilibrium product
measures.''
Let $M_1$ be the
space of functions $\gamma\dvtx \R\rightarrow[0,1]$, and let
$M_1(\rho_*,\rho^*)$ be those functions in $M_1$ which equal $\rho_*$
for all $x\leq x_*$, and which equal $\rho^*$ for all $x\geq x^*$, for
some $x_*\leq x^*$.

We will consider on $M_1$
the topology induced by $C_K(\R)$, the set of continuous, compactly
supported functions on $\R$, with the duality $\langle\cdot; \cdot
\rangle$ where $\langle\gamma; G\rangle= \int G(x)\gamma(x)\,dx$ for
$\gamma\in M_1$ and $G\in C_K(\R)$. This topology, if $M_1$ is
thought of as a measure space, is the vague topology which is metrizable.

\textit{Local equilibrium measure (LEM).}
For $0<\rho_*,\rho^* < 1$, let $\gamma\in M_1(\rho_*,\rho^*)$
be a piecewise continuous function, such that
$0<\gamma(x)<1$ for all $x\in\R$.
With respect to\vadjust{\goodbreak} $\gamma$ and a scaling parameter $N\geq1$, we define
a sequence of local equilibrium product measures $\nu^{(N)}_{\gamma
(\cdot)}$
as those formed from the marginals
$\nu^{(N)}_{\gamma(\cdot)}(\eta(x)=1) = \gamma(x/N)$ for $x\neq
0$, and $\nu^{(N)}_{\gamma(\cdot)}(\eta(0)=1)=1$.

\textit{Deterministic initial configuration (DIC).}
For $0<\rho_*,\rho^*<1$, let $\gamma$
be a piecewise continuous function in $M_1(\rho_*,\rho^*)$.
Then, a sequence of deterministic initial configurations $\xi^{\gamma,
N}$ is one such that $\xi^{\gamma, N}(0)=1$ and for all
continuous, compactly supported $G$,
$
\lim_{N\rightarrow\infty}
\frac{1}{N}\sum_{x} \xi^{\gamma, N}(x)G (x/N ) = \int
G(x)\gamma(x)\,dx.
$

We remark particular examples of local equilibrium measures
$\nu^{(N)}_{\gamma(\cdot)}$ are the equilibrium measures
$\nu_\rho(\cdot|\eta(0)=1)$ conditioned to have a particle at the
origin for $0<\rho<1$. Suitable deterministic configurations
$\xi^{\gamma, N}$, for instance, include the ``alternating'' configuration
where every other vertex is occupied corresponding to $\gamma(x)\equiv
1/2$. Nonequilibrium initial measures, corresponding to step profiles
$\gamma(x) =
\rho_*1_{(-\infty, 0]}(x) + \rho^*1_{(0,\infty)}(x)$, can also be
constructed.
The condition that the origin is occupied in these configurations
allows us to distinguish the corresponding particle as the ``tagged'' particle.

In a sense, the profiles $\gamma$, associated to the local equilibria and
deterministic profiles above, are ``nondegenerate,''
in that $\gamma$ is asymptotically
bounded strictly between $0$ and $1$.
Also, the property that $\gamma(x)$ is
constant for large $|x|$, and with respect to (LEM) specifications
that $0<\gamma<1$,
is useful to establish later
Proposition~\ref{LY_prop},
although some modifications,
for instance, in
terms of profiles sufficiently close to being constant for
large $|x|$, should be possible with more work.
However, under ``degenerate'' profiles,
different current and tagged particle large deviation behaviors might
occur. See comments after Theorems~\ref{mainthm3} and   \ref
{mainthm4} for an
``example.''

We now describe the LLNs,
proved in Jara and Landim~\cite{JL} (stated under a class of local equilibrium
measures, but the same proof also works starting from the initial
measures above):
%
\begin{equation}
\label{JL_LLN}
\lim_{N\rightarrow\infty} \frac{1}{N}J_{-1,0}(N^2t) = v_t \quad
\mbox{and}\quad\lim_{N\rightarrow\infty}\frac{1}{N}X_{N^2t}= u_t,
\end{equation}
in probability, where $v_t$ and $u_t$ satisfy
\[
\frac{dv_t}{dt} = -\frac{1}{2}\partial_x \rho(t,0) \quad\mbox
{and}\quad\frac{du_t}{dt} = -\frac{1}{2}\frac{\partial_x\rho
(t,u_t)}{\rho(t,u_t)}
\]
and $\partial_t\rho= (1/2)\partial_{xx}\rho$ and
$\rho(0,x)=\gamma(x)$, that is, $\rho(t,x) = \sigma_t*\gamma(x)$ where
$\sigma_t(y) = (2\pi
t)^{-1/2}\exp\{-y^2/2t\}$.
Note that
\[
v_t = -\frac{1}{2}\int_0^t \partial_x\rho(s,0)\,ds = \int_0^\infty
[ \rho(t,x)-\rho(0,x) ]\,dx,
\]
and $u_t$ is also the unique number $\alpha$, where
\[
\int_0^{\alpha}\rho(t,x)\,dx = -\frac{1}{2} \int_0^t \partial_x
\rho(s,0)\,ds =\int_0^\infty
[ \rho(t,x)-\rho(0,x) ]\,dx.\vadjust{\goodbreak}
\]

To explain the last equation, the right-hand side, as already
indicated, is the integrated macroscopic current across the origin up
to time $t$.
As the microscopic
dynamics is nearest-neighbor with enforced ordering of particles, the
tagged particle, initially at the origin, will be at the head of the
flow through the origin. So, to compute its macroscopic position $u_t$
at time $t$, we find $\alpha$ so that the mass at time $t$ between
positions $x=0$ and $x=\alpha$, the left-hand side of the equation,
equals the integrated current, and conclude $u_t=\alpha$.\looseness=-1

We remark, starting from a class of local equilibrium measures, corresponding
invariance principles in subdiffusive $t^{1/4}$ scale, in the sense of
finite-dimensional distributions,
with respect to fractional Brownian motion-type Gaussian processes, was
also proved in~\cite{JL}.
Also, for the current, starting from a large class of product measures,
self-normalized CLTs have
been shown in Liggett~\cite{Liggett_spa} and Vandenberg-Rodes
\cite{Vandenberg}.

In this context, we derive large deviation principles (LDPs)
(Theorem~\ref{mainthm1}), in diffusive scale, corresponding to the
laws of large numbers (\ref{JL_LLN}) when starting from (LEM) or
(DIC) measures. We give also lower and upper bounds on the
associated rate functions, starting from various nondegenerate
initial conditions (Theorem~\ref{mainthm2}).
A consequence of these
rate function bounds, say when starting from deterministic initial
configurations, is that the following growth structure can be
deduced: Namely, the rate functions are quadratic near their zeroes,
but are third order far away from the zeroes.

In particular, the third order asymptotics we derive confirm the
formal third-order expansions in Derrida--Gerschenfeld
\cite{Derrida} for the probability distribution of the current
across the origin at large times; cf. discussion after Theorem
\ref{mainthm2}. On the other hand, starting from a ``degenerate''
deterministic initial configuration with $\gamma(x) =
1_{[-1,1]}(x)$, we show that the large deviations behavior is, at
most, quadratic (Theorem~\ref{mainthm4}).

Moreover, in Theorem~\ref{mainthm3}, starting under deterministic
configurations when $\gamma(x) \equiv\rho$, we find the exact
asymptotic behavior of the rate functions near their zeroes.

The main idea for the LDPs
is to relate, through several ``entropy'' and ``energy'' estimates, the
current and tagged particle deviations to those
established in Kipnis, Olla and Varadhan~\cite{KOV}, Landim \cite
{Lan} and
Landim and Yau~\cite{LY}, with respect to the hydrodynamic limit of
the process empirical
density; cf. Propositions~\ref{hyd_limit},~\ref{KOV}. The growth
order asymptotics are proved in part by estimations
of currents and calculus of variations arguments.

At this point, we remark that the behavior of the tagged particle, in
contrast to the subdiffusive $d=1$ nearest-neighbor result, scales
differently in symmetric
exclusion models in $d\geq2$, and also in $d=1$ when the underlying
jump probability is not nearest-neighbor, that is, when particles are
free to pass by other particles. Namely, in Kipnis and Varadhan \cite
{KV}, starting under an equilibrium $\nu_\rho\{\cdot|\eta(0)=1)$,
in diffusive scale, invariance principles for the tagged particle to
Brownian motion were proved. Later, in Rezakhanlou~\cite{Rezakhanlou},
starting from local equilibrium measures,\vadjust{\goodbreak} in diffusive scale, an
invariance principle with respect to a diffusion with a drift given in
terms of the profile $\gamma$ is proved for the
``averaged'' tagged particle position, averaging over all the positions
of $O(N)$ particles in a sequence of tori with $N$ vertices. In
Quastel, Rezakhanlou and Varadhan~\cite{QRV}, in $d\geq3$, a
corresponding large deviations principle is proved for the
``averaged'' tagged particle
position with rate function, which is finite on
processes with finite relative entropy, with respect to diffusions
which, in some sense, add an additional drift to the limit diffusion in
\cite{Rezakhanlou}.
This LDP for the ``averaged'' tagged particle would seem also to hold
in $d\leq2$ (nonnearest-neighbor in $d=1$), given regularity results
on the self-diffusion coefficient in Landim, Olla and Varadhan~\cite{LOV}
not available when~\cite{QRV} was written.

We also mention, other large deviation
works with respect to empirical densities and currents in related interacting
systems are
Benois, Landim and Kipnis~\cite{BLK},
Bertini et al.~\cite{BDGJL,BDGJL1},
Bertini, Landim and Mourragui~\cite{BLM},
Farfan, Landim and Mourragui~\cite{FLM}, Quastel~\cite{Quastel}, and
Grigorescu~\cite{Grigorescu}; see also Kipnis and Landim~\cite{KL}, Chapter 10, anderences therein.
Also, we note, with respect to totally
asymmetric nearest-neighbor exclusion in $d=1$, large
deviation ``lower tail'' bounds for tagged particles are found in
Sepp\"al\"ainen~\cite{Sepplargedeviations}.

We now give the hydrodynamic limit and rate function for
the process empirical density $\mu^N(s,x;\eta)\in D([0,T]; M_1)$,
\[
\mu^N(s,x;\eta) = \sum_{k\in
\Z}\eta_{N^2s}(k)1_{[k/N,(k+1)/N)}(x)
\]
where $x\in\R, s\in[0,T]$,
and $0<T<\infty$ is a fixed time.

\begin{proposition}
\label{hyd_limit}
Starting from local equilibrium measures or deterministic
configurations,
we have for $t\in[0,T]$, $\epsilon>0$ and smooth, compactly
supported $\phi$, that
\[
\lim_{N\uparrow\infty}P \biggl\{ \biggl|\int\phi(x)\mu^N(t,x)\,dx - \int
\phi(x)m(t,x)\,dx \biggr|>\epsilon\biggr\} = 0,
\]
where $m$ satisfies $\partial_t m =
(1/2)\partial_{xx}m$ with initial data $m(0,x) =\gamma(x)$.
\end{proposition}
Aerence for the proof of Proposition~\ref{hyd_limit}, among other
places,
is Theorem~8.1 in Sepp\"al\"ainen~\cite{Sepp_book}.

The rate functions for the process empirical density differ depending
on the type of initial distribution. First, following~\cite{KOV,Lan}, suppose the process
starts from a local equilibrium measure
$\nu^{(N)}_{\gamma(\cdot)}$. For
$\mu\in D([0,T]; M_1)$, define the
linear functional on $C_K^{1,2}([0,T]\times\R)$:
\begin{eqnarray*}
l(\mu;G) &=& \int G(T,x)\mu_T(x)\,dx - \int G(0,x)\mu_0(x)\,dx \\
&&{}-
\int_0^T\int\mu_t(x) \biggl(\frac{\partial}{\partial t}
+\frac{1}{2}\,\frac{\partial^2}{\partial x^2} \biggr)G(t,x)\,dx\,dt.
\end{eqnarray*}
Let
\begin{eqnarray*}
I_0(\mu)&=&\sup_{G\in
C_K^{1,2}([0,T]\times\R)}\biggl \{l(\mu;G)-\frac{1}{2}\int_0^T\int\mu
_t(1-\mu_t)(x)G^2_x(t,x)\,dx\,dt \biggr\},\\
h(\mu_0;\gamma)&=&\sup_{\phi_0,\phi_1\in
C_K(\R)} \biggl\{\int\mu_0(x)\phi_0(x)\,dx +\int\bigl(1-\mu_0(x)\bigr)\phi_1(x)\,dx
\\
&&\hspace*{55pt}{} - \int\log\bigl[\gamma(x)
e^{\phi_0(x)} + \bigl(1-\gamma(x)\bigr)e^{\phi_1(x)}\bigr]\,dx \biggr\},
\end{eqnarray*}
and form the rate function
\[
I^{\mathrm{LE}}_{\gamma}(\mu) = I_0(\mu) + h(\mu_0;\gamma).
\]
Here, $C^{\alpha,\beta}_K$ is the space of compactly supported functions,
$\alpha$ and $\beta$-times continuously differentiable in $t$ and $x$,
respectively. In addition, we will use the notation $\mu_t(x) = \mu(t,x)$.

Next, starting from deterministic configurations
$\xi^{\gamma, N}$,
the rate
function in~\cite{LY} (written for zero-range systems, but the methods
straightforwardly apply to our exclusion context) is given by
\[
I^{\mathrm{DC}}_\gamma(\mu) = \cases{
I_0(\mu),&\quad when $ \mu_0=\gamma$,\cr
\infty,&\quad otherwise.}
\]
To simplify notation, we call both $I^{LE}_\gamma$ and
$I^{\mathrm{DC}}_\gamma$ as $I_\gamma$, omitting the super
scripts ``LE'' and
``DC,''
when statements apply to both and the context clear. For $0\leq
\alpha,\beta\leq1$, let $h_d(\alpha;\beta)=\alpha\log[\alpha
/\beta]
+ (1-\alpha)\log[(1-\alpha)/(1-\beta)]$ with usual conventions
$0\log0=0/0=0$ and $\log0 = -\infty$.

From the definition, $I_\gamma$ is a convex function. Also, a main
point in~\cite{KOV} was to note that when
$I_\gamma(\mu)<\infty$ is finite, that first
\[
h(\mu_0;\gamma) = \int h_d(\mu_0(x);\gamma(x))\,dx <\infty.
\]
[Of course, starting from deterministic
configurations, $\mu_0= \gamma$.]
Also second, $\mu$ corresponds to a function $H_x\in
L^2([0,T]\times\R,\mu(1-\mu)\,dx\,dt)$ and satisfies a ``weakly asymmetric
hydrodynamic equation,''
%
\begin{equation}\label{weak_asym_eq}\partial_t \mu= \tfrac
{1}{2}\partial_{xx}\mu-
\partial_x [H_x \mu(1-\mu) ]
\end{equation}
in the weak sense. That is, for $G\in C^{1,2}_K([0,T]\times\R)$, we have
%
\begin{equation}
\label{weak_eq}
l(\mu;G) = \int_0^T\int G_xH_x\mu(1-\mu)(t,x)\,dx\,dt
\end{equation}
and
%
\begin{equation}
\label{I_0}
I_0(\mu) = \frac{1}{2}\int_0^T\int H_x^2\mu(1-\mu)\,dx\,dt.
\end{equation}
Reciprocally, if for a density $\mu\in D([0,T]; M_1)$, there exists
$H_x\in L^2([0,T]\times\R, \mu(1-\mu)\,dx\,dt)$, such
that $\mu$ satisfies (\ref{weak_asym_eq}) weakly, then $I_0(\mu)$ is
given by~(\ref{I_0}).

Recall a function $\mathcal I\dvtx \mathcal X\to[0,\infty]$ on a complete,
separable metric space $\mathcal X$ is a rate function if it has closed
level sets $\{x\dvtx  \mathcal I(x)\leq a\}$. It is a \textit{good} rate
function if the level sets are also compact. Also,
a sequence $\{X_n\}$ of random variables with values in $\mathcal X$
satisfies a large deviation principle
(LDP) with speed $n$ and rate function $\mathcal I$
if for every Borel
set $U\in{\mathcal B}_{\mathcal X}$,
\begin{eqnarray*}
-
\inf_{x\in\bar U} \mathcal I(x) &\geq& \limsup_{n\to\infty}\frac
1n\log\Pr( X_n\in U)\\[-2pt]
&\geq& \liminf_{n\to\infty}\frac1n\log\Pr( X_n\in
U) \geq - \inf_{x\in
U^\circ} \mathcal I(x),
\end{eqnarray*}
where $U^\circ$ is the interior of $U$, and $\bar{U}$ is the closure
of $U$.

Let $\mathcal A= \mathcal A(\gamma)$ be the space of all densities
$\mu$, such that $I_\gamma(\mu)<\infty$, which can be approximated in
$D([0,T]; M_1)$ by a sequence
of densities $\{\mu^n\}$ satisfying (\ref{weak_asym_eq})
corresponding to $\{H_x^n\}\subset
C^{1,2}_K([0,T]\times\R)$, such that
$I_\gamma(\mu^n)\rightarrow I_\gamma(\mu)$.

For general local equilibrium measures (LEM) and deterministic initial
configurations (DIC),
only a weak large deviation principle is available. The next
proposition follows straightforwardly from the methods of~\cite{KOV}
(see also~\cite{KL}, Chapter 10), and replacement estimates in
\cite{LY}, namely Theorem~6.1 and Claims 1, 2~\cite{LY}, Section
6.\vspace*{-3pt}
%
\begin{proposition}
\label{LY_ldp}
With respect to initial local equilibrium measures (LEM) or
deterministic configurations (DIC), corresponding to profile $\gamma$,
$I_\gamma$ is a good convex rate
function,
and for $U\subset D([0,T]; M_1)$,
\begin{eqnarray*}
- \inf_{\mu\in\bar{U}}I_\gamma(\mu) &\geq&
\limsup_{N\uparrow\infty}\frac{1}{N}\log P [\mu^N\in
U ]\\[-2pt]
&\geq& \liminf_{N\uparrow\infty}\frac{1}{N}\log P [\mu^N\in
U ] \geq-\inf_{\mu\in U^\circ\cap\mathcal A}I_\gamma(\mu).\vspace*{-3pt}
\end{eqnarray*}
\end{proposition}

The last proposition raises the question when $\mathcal A(\gamma)$ is
large enough so that the lower bound matches the upper bound. However,
with respect to the profiles considered, the following containment is
true, so that, as a corollary, the full LDP holds.\vspace*{-3pt}
%
\begin{proposition}
\label{LY_prop}
With respect to profiles $\gamma$ associated to
local equilibrium measures (LEM) and deterministic
configurations (DIC),
\[
\mathcal A(\gamma) \supset\{\mu\dvtx  I_\gamma(\mu)<\infty
\}.\vspace*{-3pt}
\]
\end{proposition}

\begin{corollary}
\label{KOV}
With respect to initial local equilibrium measures (LEM) and
deterministic configurations (DIC), the LDP
with speed $N$ holds for $\{\mu^N\}$ with good convex rate function
$I_\gamma$.\vadjust{\goodbreak}
\end{corollary}

We note Proposition~\ref{LY_prop}, for continuous
profiles $\gamma\in M_1(\rho,\rho)$ with \mbox{$0<\rho<1$} and
$0<\gamma(\cdot)<1$ corresponding to local equilibrium measures, was
stated in~\cite{Lan}, and the associated LDP in Corollary~\ref{KOV}
with respect to these initial measures is Theorems 3.2, 3.3
\cite{Lan}. In Section~\ref{appendix}, we prove Proposition
\ref{LY_prop}, generalizing the initial states allowed.

It will be convenient to rewrite (\ref{weak_asym_eq}) in terms of
a macroscopic ``current'' or ``flux'' $J$: That is, when $I_\gamma(\mu
)<\infty$, define $J$ so that weakly,
\[
\partial_x J + \partial_t \mu= 0; \qquad  J = -\tfrac{1}{2}\partial_x \mu+
H_x\mu(1-\mu).
\]
It turns out such currents have nice properties and relations; cf.
Propositions~\ref{energyLemma} and~\ref{Lip-J}. Namely, the time
integrated current $x \mapsto\int_0^T J(x,t)\,dt$ is a well-defined
function on $\R$.
Also, the
limit
%
\begin{equation}
\label{converges_eqn}
\int_0^\infty [\mu_T(x) - \mu_0(x) ]\,dx :=\lim_{L\rightarrow
\infty}
\int_0^L [\mu_T(x) - \mu_0(x) ]\,dx \mbox{ converges}
\end{equation}
and
%
\begin{equation}\label{current-mass}
\int_0^T J(0,t)\,dt = \int_0^\infty\mu_T(x) - \mu_0(x)\,dx.
\end{equation}
In addition, for $\alpha,\beta\in\R$,
$\int_0^T [J(\beta,t)-J(\alpha,t) ]\,dt = \int_\beta^\alpha [\mu
_T(x)-\mu_0(x) ]\,dx$.

We now write the current and tagged particle rate function in terms of~$I_\gamma$.
Define
the functions $\J=\J_\gamma$ and $\I=\I_\gamma$, for $a\in\R$, by
\begin{eqnarray*}
\J(a) &=& \inf\biggl\{I_\gamma(\mu)\dvtx  \int_0^T J(0,t)\,dt = a \biggr\}\\
&=& \inf\biggl\{I_\gamma(\mu)\dvtx  \int_0^\infty\mu_T(x) - \mu_0(x)\,dx = a
\biggr\}
\end{eqnarray*}
and
\begin{eqnarray*}
\I(a) & = & \inf\biggl\{I_\gamma(\mu)\dvtx  \int_0^T J(0,t)\,dt =
\int_0^{a}\mu_T(x)\,dx \biggr\}\\
&=&
\inf\biggl\{I_\gamma(\mu)\dvtx  \int_0^\infty\mu_T(x) - \mu_0(x)\,dx =
\int_0^{a}\mu_T(x)\,dx \biggr\}.
\end{eqnarray*}
When starting from (LEM) or (DIC) initial conditions, we sometimes
distinguish the corresponding rate functions by adding a superscript.

It follows from the definitions that
%
\begin{equation}
\label{initial_ineq}
\I^{LE}_\gamma(a) \leq\I^{DC}_\gamma(a) \quad\mbox{and}\quad\J
^{LE}_\gamma(a) \leq\J^{DC}_\gamma(a).
\end{equation}
We also observe that the restriction in the infimum in the definition
of $\I$ may take different form. For instance,
when $\int_0^T J(0,t)\,dt = \int_0^a \mu_T(x)\,dx$, by the relation
$
\int_0^T J(0,t)-J(a,t)\,dt = \int_0^a \mu_T(x)-\mu_0(x)\,dx$,
one obtains the following restriction which could be used instead:
$\int_0^TJ(a,t)\,dt = \int_0^a\mu_0(x)\,dx$.\vadjust{\goodbreak}

In addition, by
translation-invariance, considering $\mu'(t,x) = \mu(t,x+a)$,
$J'(x,t)= J(x+a,t)$ and $\gamma'(x) = \gamma(x+a)$, we see,
starting from a (DIC) initial state, that
%
\begin{eqnarray}
\label{DIC_relation}
\I^{DC}_\gamma(a) & =& \inf\biggl\{I^{DC}_{\gamma}(\mu)\dvtx  \int
_0^TJ(a,t)\,dt = \int_{0}^{a}\gamma(x)\,dx \biggr\}
\nonumber\\[-9pt]\\[-9pt]
& = & \J^{DC}_{\gamma'}\biggl(\int_0^a\gamma(x)\,dx\biggr).\nonumber
\end{eqnarray}

Although one can readily see $\J$ is convex, given $I_\gamma$ is
convex and the constraint in the
definition of $\J$ is linear in $\mu$ and $a$, it is not so easily
seen whether
$\I$ is convex from this sort of argument.
However, as
seen later in Theorems~\ref{mainthm2} and~\ref{mainthm3},
near their zeroes, both $\J$ and $\I$ behave quadratically.

Also, it is perhaps curious to note that $\J$ and $\I$ can be
written
completely in terms of densities $\mu$, a consequence of the enforced
ordering of particles in the nearest-neighbor $d=1$ setting. In
contrast,
the large
deviation rate function for the ``averaged'' tagged particle position
in~\cite{QRV}
involves an auxiliary current in its description.

We now give some properties of $\J$ and $\I$ and state the large
deviation principles.

\begin{theorem}
\label{mainthm1} With respect to (DIC) or (LEM) initial measures:
\begin{longlist}[(ii)]
\item[(i)] $\J$ and $\I$ are finite on $\R$, $\lim_{|a|\uparrow\infty}
\J(a)=\lim_{|a|\uparrow\infty}\I(a) =\infty$, and $\J$ and $\I$
are a good rate functions. Further, $\J$ and $\I$ have unique zeroes
at the LLN constants $v_T$ and $u_T$, respectively.

\item[(ii)] The scaled quantities $\{J_{-1,0}(N^2T)/N\}$ and $\{X_{N^2T}/N\}$
satisfy LDPs in scale $N$ with respective rate functions $\J$ and $\I$.
\end{longlist}
\end{theorem}

A natural question at this point is to calculate the rate functions
$\J$ and $\I$.
Although this appears difficult, some bounds (with nonoptimal
constants) are possible under various conditions.

\begin{theorem}
\label{mainthm2}
Starting under (DIC) or (LEM) initial conditions, there is a constant
$c_1=c_1(\gamma)$, such that
\begin{eqnarray*}
\limsup_{a\rightarrow v_T} \frac{\sqrt{T}}{(a-v_T)^2}\J(a), \qquad
\limsup_{|a|\uparrow\infty} \frac{T}{|a|^3}\J(a) &\leq& c_1,\\[-2pt]
\limsup_{a\rightarrow u_T} \frac{\sqrt{T}}{(a-u_T)^2}\I(a),\qquad
\limsup_{|a|\uparrow\infty} \frac{T}{|a|^3}\I(a) &\leq& c_1.
\end{eqnarray*}

Also, starting under (DIC) initial conditions, there is a constant
$c_2=c_2(\gamma)>0$, such that
\begin{eqnarray*}
\liminf_{a\rightarrow v_T} \frac{\sqrt{T}}{(a-v_T)^2}\J(a), \qquad
\liminf_{|a|\uparrow\infty}\frac{T}{|a|^3}\J(a) &\geq& c_2, \\[-2pt]
\liminf_{a\rightarrow u_T}\frac{\sqrt{T}}{(a-u_T)^2}\I(a), \qquad
\liminf_{|a|\uparrow\infty}\frac{T}{|a|^3}\I(a) &\geq& c_2.\vadjust{\goodbreak}
\end{eqnarray*}
\end{theorem}

We remark the quadratic asymptotics for $\J(a)$ and $\I(a)$ near
their zeroes recalls Gaussian expansions, and the CLTs in~\cite{JL},
\cite{Liggett_spa} and~\cite{Vandenberg}. On the other hand, the
cubic bounds for large $|a|$ in Theorem~\ref{mainthm2} seem
intriguing, perhaps connected with totally asymmetric nearest-neighbor
exclusion (TASEP) effects. That is, for the current or tagged particle
to deviate to a far level $aN$, order $O(|a|N)$ particles must be
driven far away from their initial positions, so that perhaps the
process behaves like a driven system like TASEP.

We remark on these last points
that
in Derrida and Gerschenfeld~\cite{Derrida,Derrida2}, starting
from a local
equilibrium measure with step profile $\gamma^{\rho_1,\rho_2}(x) =
\rho_l1_{(-\infty,0]} + \rho_r1_{(0,\infty)}$, the large deviation
``pressure'' of the current $J_{0,1}(t)$ across the bond $(0,1)$,
$\lim_{t\uparrow\infty}t^{-1/2}\log E[\exp\{\lambda J_{0,1}(t)\}] =
F(\rho_l,\rho_r,\lambda)$, is found.
Also,
formal
asymptotics with $F$ give
$P(J_{0,1}(t) = a) \sim\exp[\sqrt{t}\{-\frac{\pi^2}{12}a^3 +
\cdots\}]$, for large $t$ and large $a>0$ (cf. page 980~\cite{Derrida2}).

In this context,
the large deviation principle in Theorem~\ref{mainthm1} and bounds in
Theorem~\ref{mainthm2} prove the form of this expression with respect
to the dominant third order term when starting from (DIC) initial
conditions: Namely, for large $a$ and constants $c_0,c_1$,
\begin{eqnarray*}
-c_0 |a|^3 & \geq& -\inf_{|x|\geq a}\J(x)
\geq\limsup_{t\rightarrow\infty}
\frac{1}{\sqrt{t}}P\bigl(|J_{0,1}(t)|\geq a\bigr)\\
&\geq&
\liminf_{t\rightarrow\infty} \frac{1}{\sqrt{t}}\log P\bigl(|J_{0,1}(t)|
\geq a\bigr) \geq-\inf_{|x|> a}\J(x) \geq-c_1 |a|^3.
\end{eqnarray*}
This addresses, in part, a question in~\cite{Derrida}, as to whether
the large $|a|$ asymptotics would extend to nonstep profiles. See also
Hurtado and Garrido~\cite{Garrido}.

Also, with respect to the current and tagged particle,
fluctuations in the ``KPZ'' class are discussed in Praehofer and Spohn
\cite{PS}, Ferrari and Spohn~\cite{FS} and Sasamoto~\cite{Sasamoto},
with respect to TASEP starting initial conditions with step or constant
profiles. In particular, the scaling
limits of the current and tagged particle are of ``Tracy--Widom'' or
``Airy'' process types whose marginal distribution have upper tail
on
order $e^{-c_0|x|^3}$
as $x\uparrow\infty$, and lower tail on order
$e^{-c_1|x|^{3/2}}$ as $x\downarrow\infty$, for some constants $c_0,c_1$.
In our context, starting from (DIC) initial conditions, we have from
Theorem~\ref{mainthm2} that $\J(a), \I(a)$ are on cubic order
$|a|^3$ for large $|a|$. Formally, one is tempted to link this cubic order
in terms of the TASEP scaling limit
process exponents.
It would be interesting
to investigate such analogies.

We nowine the behavior of $\J(a)$ and $\I(a)$ near their zeroes
$v_T=u_T=0$ when the deterministic initial condition has constant
profile $\gamma\equiv\rho$.
Arratia's CLT variances $\sigma^2_J$ and $\sigma_X^2$, mentioned
earlier, can be computed by
adding static and dynamic contributions, due to initial configuration
and later motion fluctuations, respectively.
However, starting
from deterministic initial configurations, only the dynamical
contributions would be
present, and we show later, in Proposition~\ref{dynamical_part}, that these
parts of the variances are $\sigma_{J,dyn}^2= \sqrt{\pi}\rho(1-\rho
)$ and $\sigma_{X,dyn}^2 = \sqrt{\pi}(1-\rho)/\rho$.

\begin{theorem}
\label{mainthm3}
For $\rho\in(0,1)$, starting from (DIC) initial configurations with
profile $\gamma\equiv\rho$, we have
\[
\lim_{|a|\downarrow0}\frac{1}{a^2}\J(a) = \frac{1}{2\sigma
_{J,dyn}^2\sqrt{T}} = \frac{\sqrt{\pi}}{2\sqrt{T}}\rho(1-\rho)
\]
and
\[
\lim_{|a|\downarrow0}\frac{1}{a^2}\I(a) = \frac{1}{2\sigma
_{X,dyn}^2\sqrt{T}} = \frac{\sqrt{\pi}}{2\sqrt{T}}\frac{\rho
}{1-\rho}.
\]
\end{theorem}

At this point, one might ask about the large deviation behavior
starting from initial conditions with ``degenerate'' profiles.
In this case, diffusive scaling may not
always capture for the tagged particle nontrivial LLNs, as in (\ref
{JL_LLN}) or large
deviations as in Theorem~\ref{mainthm1}. For instance, starting under
$\xi^{\gamma,N}$ where
$\gamma(x)=1_{(-\infty,0]}(x)$ is the step profile, in Arratia
\cite{Arratia} it is shown that $t^{-1/2}x(t) - \sqrt{\log(t)}
\rightarrow0$ a.s. which shows that the tagged particle diverges at
rate $\sqrt{t\log(t)}$. With respect to large deviations, it is clear
the tagged particle, initially at the origin, cannot travel to
negative locations. Also, for $a\geq0$, the condition in $\I(a)$
reduces to $\int_a^\infty\mu_T(x)\,dx = 0$ which, given that $\mu(t,x)$
satisfies (\ref{weak_asym_eq}), is impossible since the density
formally becomes positive on $\R$ as soon as $t>0$. Hence, starting
from this
step profile configuration, formally $\I= \infty$. However, for the
current, starting from this initial condition, in diffusive scaling,
$v_T<\infty$, and a corresponding CLT is proved in~\cite{Liggett_spa}.

On the other hand, when the degenerate initial profile has a density of
particles around the tagged particle, diffusive scaling would still
seem appropriate to establish an LDP for the tagged particle and
current. Here, as a contrast to the results in Theorem~\ref{mainthm2}
and to argue this last sentiment, we
show quadratic upper bounds for the current and tagged particle large deviations
starting from the degenerate
configuration $\xi^{\gamma_1, N}$ where $\xi^{\gamma_1,N}(x) = 1$ for
$|x|\leq N$ and $\xi^{\gamma_1,N}(x)=0$ otherwise. Here, $\gamma
_1(x) =
1_{[-1,1]}(x)$. Note the associated LLN speeds $v_T=u_T=0$.

\begin{theorem}
\label{mainthm4}
Starting under $\xi^{\gamma_1,N}$, there exists $c_1=c_1(T)>0$ such
that, for
$a\geq0$,
\begin{eqnarray*}
\limsup_{N\uparrow\infty} \frac{1}{N}\log P\bigl(
|J_{-1,0}(N^2T)|/N \geq a\bigr) & \leq&
\cases{ -c_1a^2, &\quad for $0\leq a \leq1$,\cr
-\infty, &\quad for $a>1$,}
\\
\limsup_{N\uparrow\infty} \frac{1}{N}\log P\bigl( |X(N^2 T)|/N
\geq
a\bigr)
& \leq& -c_1a^2.
\end{eqnarray*}
\end{theorem}

The interpretation, for instance, with respect to the tagged particle,
is that in
configurations $\xi^{\gamma_1, N}$, although it
is trapped in the middle of a large
segment of particles, to displace large distances, as there are only
$O(N)$ number of particles in the system, the cost is not as
great as under $\xi^{\rho, N}$, where there are an infinite number of
particles. At the same time, there is a positive density of particles
to the
left and right of the origin, unlike for the profile
$\gamma(x)=1_{(-\infty,0]}(x)$, which slows down the tagged particle
so that deviations to $a\in\R$ have finite cost in diffusive scale.
With respect to the current, a similar explanation applies; we note,
however, current levels larger than $N$ cannot happen, and so they are
given infinite cost.

Finally, we remark on some natural questions.

(1) As indicated by
Theorem~\ref{mainthm4}, different large deviation behaviors might
arise when starting from degenerate initial conditions. It would be of
interest to investigate these phenomena and provide estimates for the
corresponding rate functions. When starting from a degenerate initial
profile, with a density of mass around the initial tagged particle
position, although the basic argument of Theorem~\ref{mainthm1}(ii)
in Section~\ref{proofofTheorem1} holds, main obstacles are to extend
approximation Propositions~\ref{LY_prop} and~\ref{claudio}, energy
estimate Proposition~\ref{energyLemma}, first bounds and development
of the rate functions in Section~\ref{finiteness} and exponential
tightness Lemma~\ref{exponentialtightnessLemma}.

(2) Also, a joint large deviations principle for the current and tagged
particle, with rate
\[
\mathbb K(a,b) := \inf\biggl\{ I_\gamma(\mu) \dvtx  \int_0^\infty\bigl (\mu
_T(x) - \mu_0(x) \bigr)\,dx = \int_0^b \mu_T(x)\,dx = a\biggr\},
\]
should hold by the methods of the article. In this case, asymptotics of
the rate function $\mathbb K(a,b)$ for $(a,b)$ near $(v_T,u_T)$ might
be studied.

The plan of the paper is now to develop preliminary estimates in
Section~\ref{firstestimatessection}. In Section
\ref{proofofTheorem1}, we prove
Theorem~\ref{mainthm1}. Then, in
Section~\ref{asymptotic_eval_section}, we
prove Theorems~\ref{mainthm2},~\ref{mainthm3}
and~\ref{mainthm4}.
These last two sections
can be read independently of each other. Finally, in Section
\ref{appendix}, as remarked earlier, we prove Proposition~\ref
{LY_prop}, and other approximations.

\section{Preliminary estimates}
\label{firstestimatessection}

We develop, in several subsections, ``energy'' and current estimates
with respect to finite rate densities, and also prove that $\J$ and
$\I$ are a finite-valued rate functions.

\subsection{Approximation and limit estimates}
We state an approximation result derived in the course of the proof of
Proposition~\ref{LY_prop}, and also certain useful limits at infinity.
Proofs of these results are given in Section~\ref{appendix}.

\begin{proposition}
\label{claudio} Let $\mu$ be a density such that $I_0(\mu)<\infty$.
Then for all $\epsilon>0$,
there is $\mu^+\in D([0,T]; M_1)$, such that:
\begin{longlist}[(viii)]
\item[(i)] $\exists0<\delta<1$ such that $\delta\leq\mu
^+(t,x)\leq1-\delta$ for $(t,x)\in[0,T]\times\R$,\vadjust{\goodbreak}
\item[(ii)] $\mu^+\in C^{\infty}([0,T]\times\R)$,
\item[(iii)] $H_x^+\in C_K^{\infty}([0,T]\times\R)$ and
\item[(iv)]
$\|\partial^{(k)}_x\,\partial_t^{(l)}\mu^+\|_{L^\infty([0,T]\times
\R)}<\infty$ for $k,l\geq1$.
\item[(v)] If $\mu_0\equiv\gamma\in M_1(\rho_*,\rho^*)$, then
$\mu_0^+ = \sigma_{\alpha}*\gamma$ for an $\alpha>0$. In
particular, if $\mu_0(x)\equiv\rho$, then $\mu_0^+(x)\equiv
\rho$.
\item[(vi)] Also, Skorohod distance $d(\mu^+, \mu)<\epsilon$ in
$D([0,T]; M_1)$,
\item[(vii)]
$|I_0(\mu^+)- I_0(\mu)|<\epsilon$.
\item[(viii)] Also, suppose
$\hat{\gamma} \in M_1(\rho_*,\rho^*)$ is piecewise continuous, and
$0<\hat{\gamma}(x)<1$ for $x\in\R$.
Then, if $h(\mu_0;\hat{\gamma})<\infty$, we have $|h(\mu^+_0;\hat
{\gamma})-
h(\mu_0;\hat{\gamma})|\leq\epsilon$.
\end{longlist}
\end{proposition}

We remark, of course, Proposition~\ref{claudio} implies that if
$I_0(\mu)<\infty$, there is a sequence of densities $\mu^n$
satisfying properties (i)--(viii) which converges to $\mu$ in
$D([0,T]; M_1)$.

\begin{lemma}
\label{cor_limit} Let $\hat\gamma\in M_1(\rho_*,\rho^*)$, and
$\mu$ be a smooth density such that $h(\mu_0;\hat\gamma)<\infty$,
$I_0(\mu)<\infty$, and which also satisfies \textup{(i)--(iv)} in Proposition~\ref{claudio}
Then, we have
\[
\lim_{|y|\uparrow\infty}\sup_{t\in[0,T]}|\mu(t,y)-\hat\gamma
(y)| = 0.
\]
\end{lemma}

The next lemma will be used in the proof of Theorem~\ref{mainthm3}.
%
\begin{lemma}\label{unif_limit} Let $\{\mu\}$ be a
smooth density such that $\mu_0(x)\equiv\rho$, $I_0(\mu)<\infty$,
and which satisfies \textup{(i)--(iv)} in Proposition~\ref{claudio}.
Then
\[
\sup_{0\leq t\leq T}\int(\mu_t-\rho)^2\,dx \leq8I_0(\mu).
\]
\end{lemma}

\subsection{``Energy'' and current estimates}
\label{approx_subsection}
We give a formula for the rate $I_0(\mu)$, bounds on the ``energy'' $\|
\partial_x \mu\|_{L^2}$, and relations with the current.

\begin{proposition}
\label{energyLemma} Let $\mu$ be a smooth density, with finite rate
$I_0(\mu)$, satisfying \textup{(i)--(iv)} in Proposition~\ref{claudio}.
Suppose also there is a smooth $\hat\gamma\in M_1(\rho_*,\rho^*)$,
strictly bounded between $0$ and $1$, such that
$h(\mu_0;\hat\gamma)<\infty$.
Then,
%
\begin{eqnarray}
 \qquad I_0(\mu) &=& \frac{1}{8}\int_0^T\int
\frac{(\partial_x\mu)^2}{\mu(1-\mu)}\,dx\,dt +
\frac{1}{2}[h(\mu_T;\hat\gamma)-h(\mu_0;\hat\gamma)]
\nonumber\\
&&{} +\frac{1}{2} \int\frac{\partial_x\hat\gamma}{\hat\gamma
(1-\hat\gamma)}\int_0^T J\,dt\,dx+ \frac{1}{2}\int_0^T\int\frac
{J^2}{\mu(1-\mu)}\,dx\,dt,
\nonumber\\
\label{cor_bound}
 \qquad \frac{1}{4}\|\partial_x\mu\|^2_{L^2} &\leq& h(\mu_0;\hat\gamma)
+ I_0(\mu) + T\bigl\|\partial_x
\hat\gamma/\bigl(\hat\gamma(1-\hat\gamma)\bigr)\bigr\|_{L^2}^2
\end{eqnarray}
and
%
\begin{equation}
\label{current-massformula}
  \int_0^TJ(a,t)\,dt - \int_0^T J(b,t)\,dt = \int_a^b \mu_T(x)-\mu
_T(0)\,dx\qquad\mbox{for }a,b\in\R.\hspace*{-28pt}\vadjust{\goodbreak}
\end{equation}
\end{proposition}

\begin{pf}
First, as $J = -(1/2)\partial_x\mu+ H_x\mu(1-\mu)$, we have
\begin{eqnarray*}
I_0(\mu) &=& \frac{1}{2}\int_0^T\int H_x^2\mu(1-\mu)\,dx\,dt\\
& =& \frac{1}{8}\int_0^T\int\frac{(\partial_x\mu)^2}{\mu(1-\mu
)}\,dx\,dt\\
&&{}+\frac{1}{2}\int_0^T\int\frac{J\partial_x\mu}{\mu(1-\mu
)}\,dx\,dt + \frac{1}{2}\int_0^T\int\frac{J^2}{\mu(1-\mu)}\,dx\,dt.
\end{eqnarray*}
We now find a suitable expression for the middle term.
Let $G_L$ be a smooth, nonnegative, compactly supported function
in $[-L,L]$, bounded by $1$, which equals $1$ on $[-L+1,L-1]$, and
$\sup_L\int_{A_L}(G'_L)^2/G_L \,dx<\infty$ where $A_L =
[L-1,L]\cup[-L,-L+1]$. Then
%
\begin{eqnarray}
\label{entropy_calc}
&&\partial_t \int G_L(x) h_d(\mu_t(x); \hat{\gamma})\,dx \nonumber \\
&&\qquad= - \frac{1}{2}\int G_L(x)\frac{(\partial_x\mu_t)^2}{\mu
_t(1-\mu_t)}\,dx +
\int G_L (x) H_x \partial_x\mu_t \,dx\nonumber\\[-8pt]\\[-8pt]
&&\qquad\quad{}+ \frac{1}{2}\int G_L(x)\frac{\partial_x\hat
{\gamma}\partial_x\mu}{\hat{\gamma}(1-\hat{\gamma})}\,dx
- \int G_L(x)H_x\mu(1-\mu)\frac{\partial_x\hat{\gamma}}{\hat
{\gamma}(1-\hat{\gamma})}\,dx\nonumber\\
&&\qquad\quad{}+ \int_{A_L} G'_L(x) \bigl[-(1/2)\partial_x\mu_t +
H_x\bigl(\mu_t(1-\mu_t)\bigr) \bigr] \log\frac{\mu_t}{1-\mu_t}\frac{1-\hat
{\gamma}}{\hat{\gamma}}\,dx.\nonumber
\end{eqnarray}

Hence, by Schwarz's inequality and $0\leq\mu\leq1$, we can bound,
with respect to a universal constant $C$,
\begin{eqnarray*}
&&\int G_L(x) h_d(\mu_T(x);\hat{\gamma}(x))\,dx +
\frac{1}{4}\int_0^T\int G_L(x) \frac{(\partial_x\mu_s)^2}{\mu
_s(1-\mu_s)}\,dx\,ds
\\
&& \qquad\leq \int G_L(x) h_d(\mu_0(x);\hat{\gamma}(x))\,dx + C\int
_0^T\int
(H_x)^2\mu_s(1-\mu_s)\,dx\,ds \\
&& \qquad\quad{}+ CT\int
G_L(x)\frac{(\partial_x\hat\gamma)^2}{\hat\gamma^2(1-\hat\gamma)^2}\,dx\\
&& \qquad  \quad {}+ C\int_{A_L} [(G'_L)^2/G_L]\biggl[\log\frac{\mu_t}{1-\mu_t}\frac
{1-\hat{\gamma}}{\hat{\gamma}}\biggr]^2\,dx.
\end{eqnarray*}
We can take $L\uparrow\infty$, so that the last term
vanishes by Lemma~\ref{cor_limit}.
Then, by monotone convergence, with respect to a universal constant $C$,
\begin{eqnarray*}
&&h(\mu_T;\hat{\gamma}) +
\frac{1}{4}\int_0^T\int\frac{(\partial_x\mu_s)^2}{\mu_s(1-\mu
_s)}\,dx\,ds\\
&&\qquad\leq h(\mu_0;\hat{\gamma}) + C\int_0^T\int(H_x)^2\mu
_s(1-\mu_s)\,dx\,ds
+ CT\|\partial_x\hat\gamma\|_{L^2},
\end{eqnarray*}
and as $0<\mu, \hat\gamma<1$, we have
$\|\partial_x\mu\|_{L^2}, \|J\|_{L^2}<\infty$.

Hence, integrating (\ref{entropy_calc})
and taking limit on $L$, the middle term equals
\begin{eqnarray*}
\int_0^T\int\frac{J\partial_x\mu}{\mu(1-\mu)}\,dx\,dt &=& h(\mu
_T;\hat\gamma)-h(\mu_0;\hat\gamma)
+ \int_0^T\int\frac{J\partial_x\hat\gamma}{\hat\gamma(1-\hat
\gamma)}\,dx\,dt.
\end{eqnarray*}

The desired bound on $\|\partial_x \mu\|_{L^2}$ now follows. Since
$\mu(1-\mu)\leq1/4$ and $\|H_x\mu(1-\mu)\|_{L^2}^2 \leq I_0(\mu)$
(cf. (\ref{I_0})), by Schwarz's
inequality,
we may write
\begin{eqnarray*}
\|\partial_x\mu\|^2_{L^2} & \leq& h(\mu_0;\hat\gamma) +
\frac{1}{2}\|J\|^2_{L^2} +
\frac{T}{2}\bigl\|\partial_x\hat\gamma/\bigl(\hat\gamma(1-\hat\gamma)\bigr)\bigr\|^2_{L^2}
+ 2I_0(\mu)\\
& \leq& h(\mu_0;\hat\gamma) +
\frac{1}{4}\|\partial_x\mu\|^2_{L^2} + \frac{5}{2}I_0(\mu)+
T\bigl\|\partial_x\hat\gamma/\bigl(\hat\gamma(1-\hat\gamma)\bigr)\bigr\|^2_{L^2}.
\end{eqnarray*}

Finally, (\ref{current-massformula}) expresses that the difference of
the currents across $a$ and $b$
up to time $T$ is equal to the difference in the masses in the
interval $[a,b]$ from times $T$ to $0$. This is obtained by integrating
$\partial_x J = -\partial_t\mu$.
\end{pf}

\begin{corollary}
\label{cor_weakconvergence}
Let $\mu$ be a density with finite rate $I_\gamma(\mu)<\infty$. Let
also $\{\mu^n\}$ be a sequence converging to $\mu$ with properties
\textup{(i)--(viii)} in Proposition~\ref{claudio}. Then, $\partial_x\mu^n$ and
$J^n$ are uniformly bounded in $L^2([0,T]\times\R)$ and $\partial
_x\mu^n\rightarrow\partial_x\mu$, $J^n\rightarrow J$ weakly in
$L^2([0,T]\times\R)$; consequently, $\partial_x\mu, J\in
L^2([0,T]\times\R)$.
\end{corollary}

\begin{pf}
Let $\hat{\gamma}$ be a smooth function in $M(\rho_*,\rho^*)$ such
that $0<\gamma_*<\hat{\gamma}<\gamma^*<1$ for some constants
$\gamma_*,\gamma^*$, and $h(\gamma; \hat\gamma)<\infty$. Then, by
property (viii) Proposition~\ref{claudio}, as $h(\mu_0;\hat\gamma
)<\infty$, we have $h(\mu^n_0;\hat\gamma) \rightarrow h(\mu_0;\hat
\gamma)$, and, in particular, $\{h(\mu^n_0;\hat\gamma)\}$ is
uniformly bounded.

Also, as $I_0(\mu)<\infty$, by property (vii) Proposition \ref
{claudio}, we have $I_0(\mu^n)\rightarrow I_0(\mu)$ and $\{I_0(\mu
^n)\}$ is uniformly bounded. In particular,
$\{\|H^n_x\mu^n(1-\mu^n)\|_{L^2}\}$ is uniformly bounded.

Hence, as $\partial_x\hat\gamma/(\hat\gamma(1-\hat\gamma))\in
L^2$, and by (\ref{cor_bound}) in Proposition~\ref{energyLemma}, we have
$\{\|\partial_x\mu^n\|_{L^2}\}$ is uniformly bounded. Also, since
$J^n = (1/2)\partial_x\mu^n + H^n_x\mu^n(1-\mu^n)$, we also
conclude $\{\|J^n\|_{L^2}\}$ is uniformly bounded.

We can then extract subsequences $\partial_x\mu^{n_k}$ and $J^{n_k}$
converging weakly to $\zeta$ and $\phi$, respectively. Given $\mu
^{n_k}\rightarrow\mu$ in $D([0,T]\times M_1)$, for smooth, compactly
supported $G$, we have
$\int G\partial_x\mu^{n_k} \,dx\,ds =\int-G_x\mu^{n_k} \,dx\,ds$ converges
to both $\int G\zeta \,dx\,ds$ and $\int-G_x\mu \,dx\,ds$. Then, $\partial
_x\mu$ exists weakly in $L^2$ and $\partial_x\mu= \zeta$. Hence,
the whole sequence $\partial_x\mu^n \rightarrow\partial_x\mu$
weakly in $L^2$.

Similarly, noting Skorohod convergence $\mu^n\rightarrow\mu$ implies
at the endpoints that $\mu^n_0,\mu^n_T$ converge to $\mu_0,\mu_T$,
respectively, and $\partial_t \mu^n + \partial_x J^n =0$,
we have
$\phi_x = -\partial_t\mu$ weakly in $L^2$.
Then, $\phi_x = (-1/2)\partial_{xx} \mu+ \partial_x[H_x \mu(1-\mu
)]$ weakly in $L^2$, and so $\phi= (-1/2)\partial_x\mu+ H_x\mu
(1-\mu) + C(t)$ with respect to a function $C(t)$ not dependent on\vadjust{\goodbreak}
$x$. But, given $\phi, \partial_x\mu, H_x\mu(1-\mu)\in
L^2([0,T]\times\R)$, we conclude $C(t)\equiv0$. In particular, $\phi
= J= -(1/2)\partial_x\mu+ H_x\mu(1-\mu)\in L^2$, and the sequence
$J^n \rightarrow J$ weakly in $L^2$.\vspace*{-3pt}
\end{pf}

\subsection{Current-mass relation}

We give some properties of the integrated current $\int_0^T J(x,t)\,dt$
and prove the current-mass relation indicated in the \hyperref[sec0]{Introduction}.\vspace*{-3pt}

\begin{proposition}
\label{Lip-J}
Let $\mu$ be a density such that $I_\gamma(\mu)<\infty$.
Let $\{\mu^n\}$ be a sequence converging to $\mu$ with properties
\textup{(i)--(viii)} in Proposition~\ref{claudio}. Then, $x\mapsto\int_0^T
J(x,t)\,dt$ is a Lipschitz function, $\lim_{|x|\uparrow\infty}\int
_0^T J(x,t)\,dt = 0$, and pointwise for $x\in\R$,
\[
\lim_{n\rightarrow\infty}\int_0^T J^n(x,t)\,dt= \int_0^T J(x,t)\,dt.
\]

In addition, convergence (\ref{converges_eqn}), and the
``current-mass'' relation (\ref{current-mass})~hold.\looseness=-1\vspace*{-3pt}
\end{proposition}

\begin{pf}
First, from (\ref{current-massformula}) in Proposition \ref
{energyLemma}, we have
\[
\int_0^TJ^n(a,t)\,dt - \int_0^T J^n(b,t)\,dt = \int_a^b \mu^n_T(x)-\mu
^n_T(0)\,dx.
\]
Hence $|\int_0^TJ^n(a,t)\,dt - \int_0^T J^n(b,t)\,dt| \leq|b-a|$ as
$0\leq\mu^n\leq1$. In particular, $\int_0^T J^n(a,t)\,dt$ is
Lipschitz in $a$.
Moreover, a subsequence, $\int_0^T J^{n_k}(\cdot,t)\,dt \rightarrow
\psi(\cdot)$ converges uniformly on compact subsets to a Lipschitz
function $\psi$.
Given $J^n \rightarrow J$ weakly in $L^2([0,T]\times\R)$ by Corollary
\ref{cor_weakconvergence}, we conclude by a limit argument with
respect to $G\in L^2(\R)$ that $\int G(a)\int_0^TJ(a,t)\,dt\,da=\int
G(a)\psi(a)\,da$, and so $\psi(a) = \int_0^TJ(a,t)\,dt$. In particular,
the whole sequence $\int_0^T J^n(\cdot,t)\,dt \rightarrow\int
_0^TJ(\cdot,t)\,dt$ and the limit $\int_0^TJ(\cdot,t)\,dt$ is Lipschitz.

Therefore, since
\[
\int\biggl[\int_0^T J(x,t)\,dt \biggr]^2\,dx \leq T\int\int_0^T J^2(x,t)\,dt\,dx<\infty,
\]
we obtain the pointwise limit $\int_0^TJ(x,t)\,dt\rightarrow0$ as
$|x|\uparrow\infty$.

Finally, given Skorohod convergence $\mu^n\rightarrow\mu$, $\mu
_0^n$ and $\mu^n_T$ converge respectively to $\mu_0$ and $\mu_T$.
Then, by taking limits, we can write
\[
\int_0^T J(0,t)\,dt - \int_0^T J(L,t)\,dt = \int_0^L \mu_T(x) - \mu
_0(x)\,dx.
\]
Now, since $\lim_{L\rightarrow\infty} \int_0^T J(L,t)\,dt = 0$, we
obtain (\ref{converges_eqn}) and (\ref{current-mass}).\vspace*{-3pt}
\end{pf}

\subsection{First estimates on $\J$ and $\I$}\label{finiteness}

We develop some first bounds on $\J$ and~$\I$, and at the end show
they are rate functions.

Recall $\sigma_t(x) = (2\pi
t)^{-1/2}\exp\{-x^2/2t\}$, and consider
a $C^\infty$ smooth function, supported on $[-1,1]$, say
\[
\psi_0(x) = \exp\{-1/(1-x^2) \}.\vadjust{\goodbreak}
\]
Define the smooth, anti-symmetric function
\[
\psi(x) = \cases{
-\psi_0\bigl(2(x+1/2)\bigr),&\quad for $x\leq0$,\cr
\psi_0\bigl(2(x-1/2)\bigr),&\quad for $x\geq0$
}
\]
and also the anti-derivative $\Psi(x) = \int_{-1}^x\psi(y)\,dy$,
both supported on $[-1,1]$.

Let $\gamma\in M_1(\rho_*,\rho^*)$ be a profile associated to an
initial (LEM) local equilibrium measure or a (DIC) deterministic
configuration. Recall, when $I_\gamma(\mu)<\infty$, it has explicit
representation; cf.~\ref{I_0}.
Recall, also that $v_T$ and $u_T$ are the LLN speeds associated to
$\gamma$; cf. (\ref{JL_LLN}).

Since $\J$ and $\I$ are given through infima, it is natural to look
for explicit densities where computations can be made.
Consider the density
\[
\mu(s,x)=\sigma_s*\gamma(x) + \bigl(\lambda
\epsilon(s/T)\bigr)\psi(x/L),
\]
where $\epsilon(t)$ is a smooth, increasing function which
vanishes for $0\leq t\leq1/10$, and
$\epsilon(1)=1$, and $L\neq0$.
At time $s=T/10$, $0<\gamma_*<\sigma_{s}*\gamma<\gamma^*<1$ for
some constants $\gamma_*,\gamma^*$. We will take $0\leq\lambda<
\min\{\gamma_*, 1-\gamma^*\}/2$, small
enough so that $\gamma_*/2\leq\mu\leq(1-\gamma^*)/2$ for $T/10\leq
t\leq T$.

Then, as $\mu$ follows the heat equation for $[0,T/10]$, $\mu$
satisfies (\ref{weak_asym_eq}) with respect to
$H_x$, supported on $[T/10,T]\times[-|L|,|L|]$, given by
\[
H_x = \cases{\displaystyle
\frac{1}{\mu(1-\mu)} \biggl[ \frac{\lambda\epsilon(s/T)}{2L}\psi'
\biggl(\frac{x}{L} \biggr) -\frac{\lambda
L\epsilon'(s/T)}{T}\Psi\biggl(\frac{x}{L} \biggr) \biggr],\vspace*{4pt}\cr \quad \qquad \mbox{for
$\displaystyle\frac{T}{10}\leq s\leq T, |x|\leq|L|$},\vspace*{2pt}\cr
0,  \qquad  \mbox{otherwise}.}
\]
Also, as $\mu_0 = \gamma$, we have $h(\mu_0;\gamma)=0$,
and
%
\begin{eqnarray}\label{I_0_bound}
I_0(\mu) & = & \frac{1}{2}\int_{T/10}^T \int
\frac{1} {\mu(1-\mu)} \biggl[\frac{\lambda\epsilon(s/T)}{2L}\psi'
\biggl(\frac{x}{L} \biggr)
- \frac{\lambda L\epsilon'(s/T)}{T}\Psi\biggl(\frac{x}{L} \biggr)
\biggr]^2\,dx\,ds\nonumber\hspace*{-35pt}\\[-8pt]\\[-8pt]
& \leq&\frac{4\epsilon^*}{\gamma_*(1-\gamma^*)} \biggl[\frac{\lambda^2
T}{4 |L|}\int_{-1}^1\psi'(x)^2\,dx + \frac{\lambda^2|L|^3}{T}\int
_{-1}^1\Psi(x)^2\,dx \biggr],
\nonumber\hspace*{-35pt}
\end{eqnarray}
where $\epsilon^* = 1+\|\epsilon'\|^2_{L^\infty}$.
Compute now
\begin{eqnarray*}
\int_0^\infty [\mu_T(x)-\mu_0(x) ]\,dx&=&\lambda L\int_0^1\psi(x)\,dx
+ v_T\\
&=&\lambda L\int_0^1\psi(x)\,dx + \int_0^{u_T}\sigma_T*\gamma(x)\,dx,
\end{eqnarray*}
and, for $\c\in\R$,
\begin{eqnarray*}
\int_0^\c\mu_T(x)\,dx
&=& \int_0^\c\sigma_T*\gamma(x)\,dx + \lambda L\int_0^{|\c
|/|L|}\psi(x)\,dx.
\end{eqnarray*}

Then, the restriction specified in the definition of $\J(\c)$, $\int
_0^TJ(0,t)\,dt = \c$, is the same as
%
\begin{equation}\label{J_restriction}
\lambda L\int_0^1\psi(x)\,dx = \c - v_T,
\end{equation}
and the restriction listed in $\I(\c)$,
\[
\int_0^T J(0,t)\,dt = \int_0^\infty\mu_T(x)-\mu_0(x)\,dx = \int_0^\c
\mu_T(x)\,dx,
\]
is equivalent to
%
\begin{equation}
\label{I_restriction}
\lambda L\int_{|\c|/|L|}^1\psi(x)\,dx = \int_{u_T}^\c
\sigma_T*\gamma(x)\,dx.\vspace*{-3pt}
\end{equation}

\begin{lemma}
\label{finite_rate_prop} For $\c\in\R$,
$\J(\c), \I(\c)<\infty$ and in particular $\J(v_T)=\I(u_T)=0$.
Moreover, on any interval $[a,b]\subset\R$,
$\sup_{\c\in[a,b]}\J(c), \sup_{\c\in[a,b]}\I(c)<\infty$.\vspace*{-3pt}
\end{lemma}

\begin{pf}
For $\c\in\R$, given bound (\ref{I_0_bound}), we need only
demonstrate that restrictions (\ref{J_restriction}) and (\ref
{I_restriction}), with respect to $\J$ and $\I$, hold with respective
choices of $\lambda$ and~$L$. If $\c=v_T$ or $u_T$, we may take
$\lambda=0$, and so clearly $\J(v_T)=\I(u_T)=0$.

For
$\c\neq v_T$, let $\lambda>0$, and note the left-hand side of (\ref
{J_restriction}) can be made equal to the right-hand side $v_T-\c$
with a proper choice of $L$. Similarly, when
$\c\neq u_T$, let $\lambda>0$, and note that the left-hand side of
(\ref{I_restriction}) vanishes for $|L|\leq|\c|$ and diverges to
$\pm\infty$ as $L\rightarrow
\pm\infty$.
Hence, a proper choice of $L$ allows us to verify (\ref
{I_restriction}) also.

In particular, we can see, by varying $L$, with respect to $\c\in
[a,b]$ in any finite interval, we obtain $\sup_{\c\in[a,b]}\J(\c),
\sup_{\c\in[a,b]}\I(\c) < \infty.$\vspace*{-3pt}
\end{pf}

\begin{lemma}\label{lowersemicontinuous} With respect to local
equilibrium measures or deterministic initial configurations, $\J$ and
$\I$ are lower semi-continuous.\vspace*{-3pt}
\end{lemma}

\begin{pf} We give the proof for $\I$; the argument for $\J$ is
analogous. We first consider when starting from a local equilibrium
measure and $I_\gamma= I^{LE}_\gamma$.
Let $\{a_n\}$ be a convergent sequence $a_n
\rightarrow a$. From Proposition~\ref{finite_rate_prop}, we have $\sup
_n \I(a_n)<\infty$. Then, by
Propositions~\ref{claudio} and~\ref{Lip-J}, we can find
densities
$\{\mu^n\}$ so that $|I^{LE}_\gamma(\mu^n) - \I(a^n)|<n^{-1}$ and
\mbox{$|\int_0^TJ^n(0,t)\,dt-\int_0^{a_n}\mu^n_T(x)\,dx|\leq n^{-1}$}.

As $I^{LE}_\gamma$ is a good rate function and $\{I^{LE}_\gamma(\mu
^n)\}$ is uniformly bounded, a subsequence can be found
where $\mu^{n_k}$ converges to a density $\hat{\mu}$ in $D([0,T];
M_1)$ and $\liminf\I(a^n) = \lim\I(a^{n_k})=\lim I^{LE}_\gamma(\mu
^{n_k})$.

By Proposition~\ref{Lip-J}, we have
$\int_0^T J^{n_k}(0,t)\,dt \rightarrow \int_0^T \hat{J}(0,t)\,dt$.
Also, as $\mu^{n_k}_T \rightarrow\hat{\mu}_T$, and $a_n\rightarrow
a$, we have
$\int_0^{a_n}\mu^{n_k}_T(x)\,dx \rightarrow \int_0^a \hat{\mu}_T(x)\,dx$.
Then, $\int_0^T \hat{J}(0,t)\,dt =
\int_0^a \hat{\mu}_T(x)\,dx$, and hence $\hat\mu$ satisfies the
infimum restriction in the definition of $\I(a)$.\vadjust{\goodbreak}

By lower semi-continuity
of $I^{LE}_\gamma$, the desired lower semi-continuity of $\I$ follows
as $\liminf\I(a_n) = \lim I^{LE}_\gamma(\mu^{n_k}) \geq
I^{LE}_\gamma(\hat\mu)
\geq\I(a)$.

Starting from a deterministic configuration, we can repeat the steps
with $I^{LE}_\gamma$ replaced by $I_0$. The densities $\{\mu^{n_k}\}
$, by Proposition~\ref{claudio}, also are such that
$\mu^{n_k}_0$ converges to $\gamma$. Hence, the limit
$\hat{\mu}$ satisfies $\hat{\mu}_0 = \gamma$ and so
$I_0(\hat{\mu})=I^{DC}_\gamma(\hat\mu)$. Therefore, $\I$ is also
lower semi-continuous in this case.\vspace*{-3pt}
\end{pf}

\begin{corollary}
\label{I_ratefunction}
With respect to local equilibrium measures or deterministic initial
conditions, $\J$ and $\I$ are finite-valued rate functions. In
addition, $\J(a')=0$ and $\I(a)=0$ exactly when $a'=v_T$ and $a=u_T$.\vspace*{-3pt}
\end{corollary}

\begin{pf} We concentrate on the proof with respect to $\I$, as a
similar argument holds for $\J$.
First, that $\I\dvtx \R\rightarrow\R$, $\I(u_T)=0$, and $\I$ is a rate
function follows from Lemmas~\ref{finite_rate_prop} and \ref
{lowersemicontinuous}. We need only show that $u_T$ is the only zero of
$\I$.

When $a\neq u_T$, if $\I(a)$ vanishes, out of a
minimizing sequence of densities, through Propositions~\ref{claudio}
and~\ref{Lip-J}, one can find a subsequence converging to a minimizing
$\mu$ satisfying the restriction $\int_0^T J(0,t)\,dt =
\int_0^a \mu_T(x)\,dx$.

With respect to local equilibrium measures, by lower semi-continuity of
$h(\cdot;\gamma)$ and $I_0(\cdot)$, we have
$h(\mu_0;\gamma)=I_0(\mu)=0$. Under deterministic initial
conditions, since the subsequence at time $0$ converges to $\gamma$,
we have $\mu_0=\gamma$, and by lower semi-continuity, $I_0(\mu)=0$.

Then, in either case, $\mu_0 =
\gamma$ a.s.
and, noting (\ref{I_0}), $H_x^2\mu(1-\mu)
=0$ a.s. In particular, $\mu_t = \sigma_t*\gamma$
is the unique bounded solution of the
weak heat
equation with initial data $\gamma$.
However, then
$\int_0^TJ(0,t)\,dt =\int_0^{u_T}\mu_T(x)\,dx$ which does not equal
$\int_0^a \mu_T(x)\,dx$ since $\mu_T$ is positive and
$a\neq u_T$. This is a contradiction.\vspace*{-3pt}
\end{pf}

\section{\texorpdfstring{Proof of Theorem \protect\ref{mainthm1}}{Proof of Theorem 1.5}}
\label{proofofTheorem1}

The proofs follow in several steps which are divided
into subsections. The first step is to describe key relations between a
tagged particle and the current across the bond $(-1,0)$, which will
allow us later to invoke large deviations of the empirical density.
Next, a super-exponential inequality is given. Then, exponential
tightness is established,
and weak upper and lower large deviation
bounds are proved. Finally, Theorem~\ref{mainthm1} is shown.\vspace*{-3pt}

\subsection{Tagged particle and current relations}
\label{tagged_particle_current_subsection}
For $x\in\Z$ and $t\geq0$, define $J_{x,x+1}(t)$ as the integrated
current up to time $t$ across the bond
$(x,x+1)$, that is, the number of particles which crossed from $x$
to $x+1$ up to time $t$ minus the number of particles which moved
from $x+1$ to $x$ in time $t$.
It is well known (cf. Liggett~\cite{Liggett}, DeMasi and Ferrari
\cite{DeMF}) that for integers $r>0$,
%
\begin{equation}\label{positive_a} \{X_t \geq r \} = \Biggl\{J_{-1,0}(t)
\geq
\sum_{x=0}^{r-1}\eta_t(x) \Biggr\}.\vadjust{\goodbreak}
\end{equation}
Similarly, for $r<0$,
%
\begin{equation}
\label{negative_a}
\{X_t \leq r \} = \Biggl\{J_{-1,0}(t)
\leq-\sum_{x=r}^{-1}\eta_t(x) \Biggr\}
\end{equation}
and
\[
\{X_t \leq0 \} = \{J_{-1,0}(t) \leq0 \}.
\]
Also, from a moment's thought, we have
\[
J_{x-1,x}(N^2t) - J_{x,x+1}(N^2t) = \eta_{N^2t}(x) -
\eta_0(x).
\]

We would like to make a summation-by-parts,
\begin{eqnarray*}
J_{-1,0}(N^2t) &=& \sum_{x\geq0} J_{x-1,x}(N^2t)-J_{x,x+1}(N^2t)\\
&=& \sum_{x\geq0} \eta_{N^2t}(x) -
\eta_0(x),
\end{eqnarray*}
to write the current across the bond $(-1,0)$ in terms of the empirical
process. However, the above display is only
formal as the sum on the right may not converge. To treat it carefully,
we introduce a ``cutoff'' function as in Rost and Vares~\cite{Rost-Vares}.
For $n\geq1$, let
\[
G_n(u) = 1_{[0, n]}(u)(1-u/n).
\]
Also, denote for a function
$G\in C^\infty_K(\R)$,
\begin{eqnarray*}
Y^N_t(G) & = & \frac{1}{N}\sum_x G(x/N)\eta_{N^2t}(x).
\end{eqnarray*}
Then
\begin{eqnarray*}
Y^N_t(G_n) - Y_0^N(G_n) &=& \frac{1}{N}\sum_x
G_n(x/N) \bigl(J_{x-1,x}(N^2t)
- J_{x,x+1}(N^2t) \bigr)\\
&=&\frac{1}{N}\sum_x \bigl(G_n(x/N) - G_n(x-1/N) \bigr)J_{x-1,x}(N^2t)\\
&=& \frac{1}{N}J_{-1,0}(N^2t)
-\frac{1}{N}\sum_{x=1}^{nN}\frac{1}{nN}J_{x-1,x}(N^2t).
\end{eqnarray*}
This implies
\[
\frac{1}{N}J_{-1,0}(N^2t) = Y_t^N(G_n) - Y_0^N(G_n) +
\frac{1}{N}\sum_{x=1}^{nN}\frac{1}{nN}J_{x-1,x}(N^2t).
\]

Hence, for $a>0$,
%
\begin{eqnarray}
\label{tagged-current}
&&\{X_{N^2t}/N \geq a \}\nonumber\\
 && \qquad = \Biggl\{\frac{1}{N}J_{-1,0}(N^2t)\geq
\frac{1}{N}\sum_{x=0}^{\lfloor aN\rfloor} \eta_{N^2t}(x) \Biggr\}
\\
&& \qquad = \Biggl\{Y^N_t(G_n) - Y_0^N(G_n) + \frac{1}{nN^2}\sum
_{x=1}^{nN}J_{x-1,x}(N^2t) \geq
\frac{1}{N}\sum_{x=0}^{\lfloor aN\rfloor}\eta_{N^2t}(x) \Biggr\}
.\nonumber
\end{eqnarray}
A similar statement holds for $a\leq0$, namely,
\begin{eqnarray*} \{X_{N^2t}/N \leq a \}& = &
\Biggl\{Y^N_t(G_n) - Y_0^N(G_n) \\
&&\hphantom{\Biggl\{}{}+ \frac{1}{nN^2}\sum_{x=1}^{nN}J_{x-1,x}(N^2t) \leq
-\frac{1}{N}\sum_{x=\lfloor aN\rfloor}^{-1}\eta_{N^2t}(x) \Biggr\},
\end{eqnarray*}
where for $a=0$, we take $\sum_{x=0}^{-1}\eta_{N^2t}(x) = 0$.

Therefore, heuristically, the tagged particle large deviations
should be given in terms of the rate for the empirical density
$I_\gamma$ under a certain restriction, as long as the
contribution from the term $(1/nN^2)\sum_{x=1}^{nN}J_{x-1,x}(N^2t)$
is superexponentially small as $n,N\uparrow\infty$.

\subsection{Superexponential estimate}
In relation to (\ref{tagged-current}), the superexponential estimate
needed is implied by the following estimate.

\begin{proposition}\label{error}For each $\lambda>0$, starting from
(LEM) or (DIC) initial states,
\[
\lim_{n\uparrow\infty}\lim_{N\uparrow\infty}\frac{1}{N}\log E
{\exp}\Biggl|\frac{\lambda N}{nN^2}\sum_{x=1}^{nN}J_{x-1,x}(N^2t) \Biggr|
= 0.
\]
\end{proposition}

\begin{pf} By the inequality $e^{|x|}\leq e^x + e^{-x}$, we can
remove the absolute value in the last display. Now, note that
\begin{eqnarray*}
&&\exp\Biggl\{\frac{\lambda N}{nN^2}\sum_{x=1}^{nN}J_{x-1,x}(N^2t) -
\sum_{x=1}^{nN}
(e^{\lambda/nN}-1)\int_0^{N^2t}\eta_{x-1}(1-\eta_x)(s)\,ds\\
&&\hspace*{88pt}\hphantom{\exp\Biggl\{}{} -
\sum_{x=1}^{nN}(e^{-\lambda/nN}-1)\int_0^{N^2t}\eta_x(1-\eta_{x-1})(s)\,ds
\Biggr\}
\end{eqnarray*}
is a martingale with mean $1$.
Then together, the second and third terms in the exponent equal
\begin{eqnarray*}
&&\sum_{x=1}^{nN}\biggl[(e^{\lambda/nN}-\lambda/nN
-1)\int_0^{N^2t}\eta_{x-1}(1-\eta_x)(s)\,ds\\
&&\hphantom{\sum_{x=1}^{nN}\biggl[}{}+
(e^{-\lambda/nN}+\lambda/nN-1)\int_0^{N^2t}\eta_x(1-\eta
_{x-1})(s)\,ds\biggr]\\
&&\quad{}+ \frac{\lambda}{nN}\int_0^{N^2t}
(\eta_{0} -
\eta_{nN})(s)\,ds\\
&&\qquad\leq\frac{2e^{\lambda/nN}\lambda^2}{n^2N^2}(nN)(N^2t) +
\frac{\lambda}{nN}(N^2t)\leq\frac{C(t,\lambda)N}{n},
\end{eqnarray*}
which gives the result with standard manipulations.
\end{pf}

\subsection{Exponential tightness estimate}
\label{exponentialtightnesssection}
We now show that the scaled tagged particle positions are exponential tight.

\begin{lemma}
\label{exponentialtightnessLemma} Starting from (LEM) or (DIC) initial
states, we have
\begin{eqnarray*}
&&\lim_{a\uparrow\infty}\lim_{N\uparrow\infty} \frac{1}{N}\log P
\{|J_{-1,0}(N^2T)|/N \geq a \} \\
&&\qquad= \lim_{a\uparrow\infty}
\lim_{N\uparrow\infty} \frac{1}{N}\log P \{|X_{N^2 T}|/N \geq a \}
= -\infty.
\end{eqnarray*}
\end{lemma}

\begin{pf} We give the argument for the tagged particle, as the proof
for the current is similar, and somewhat easier.
From (\ref{tagged-current}),
we need only super-exponentially estimate, for $a$
positive (as a similar argument works for $a<0$) and $n$
fixed,
\[
P \Biggl\{ Y_T^N(G_n)-Y_0^N(G_n) +
\frac{1}{nN^2}\sum_{x=1}^{nN}J_{x-1,x}(N^2T)\geq
Y^N_T(1_{[0,a]}) \Biggr\}.
\]

We need only estimate
\begin{eqnarray*}
&&E \Biggl[\exp\Biggl\{ N\Biggl[Y_{T}^N(G_n)- Y_0^N(G_n)
+\frac{1}{nN^2}\sum_{x=1}^{nN}J_{x-1,x}(N^2T)-Y^N_T\bigl(1_{[0,a]}\bigr)\Biggr
]\Biggr\}\Biggr]\\
&&\qquad= E[ e^{Q_1}e^{Q_2}e^{Q_3}e^{Q_4}]
\end{eqnarray*}
with $Q_1 = NY_T^N(G_n)$, $Q_2 = - NY^N_0(G_n)$, $Q_3 =
(nN)^{-1}\sum_{x=1}^{nN}J_{x-1,x}(N^2t)$ and $Q_4 =
-\sum_{x=0}^{\lfloor aN\rfloor}\eta_{N^2T}(x)$.
By Chebyshev, we can estimate the exponential terms
separately. For fixed $n$, $\lim N^{-1}\log
E[e^{4Q_3}]$ is bounded from Proposition~\ref{error}, and as $Q_1 \leq
nN$ by properties of $G_n$, $\lim N^{-1}\log E[e^{4Q_1}]$ is also
bounded. In addition, as $\exp\{4Q_2\}\leq1$, this term can
be neglected.\vadjust{\goodbreak}

Finally, by Borcea, Branden and Liggett
\cite{BBL}, Theorem 5.2, as the initial measure of type (LEM) or (DIC)
is a product measure [of degenerate Bernoulli's under (DIC) initial
configurations],
the coordinates $\{\eta_{N^2 T}(x)\}$ are negatively associated. Hence,
$E[e^{4Q_4}]\leq\prod_{x=1}^{\lfloor
aN\rfloor}E[e^{-4\eta_{N^2T}(x)}]$, and using $\log(1-x)\leq-x$ for
$0\leq x\leq1$, we write
\begin{eqnarray*}
\frac{1}{N}\log E[e^{4Q_4}] &\leq& \frac{1}{N}\sum_{x=1}^{\lfloor
aN\rfloor}\log E[e^{-4\eta_{N^2T}(x)}]\\
&\leq& \frac{1}{N} \sum_{x=1}^{\lfloor
aN\rfloor}\log\bigl[ (e^{-4}-1)P\bigl(\eta_{N^2T}(x)=1\bigr)+1\bigr]\\
&\leq&\frac{e^{-4}-1}{N}E\Biggl[\sum_{x=1}^{\lfloor
aN\rfloor}\eta_{N^2T}(x)\Biggr] \rightarrow(e^{-4}-1)\int_0^a m(T,x)\,dx,
\end{eqnarray*}
where $m(T,x) = \sigma_T*\gamma(x)$ is the solution of the
hydrodynamic equation (Proposition~\ref{hyd_limit}). Since $\sigma
_{T/10}*\gamma(x) \geq
\gamma_*>0$ as $\gamma\in M_1(\rho_*,\rho^*)$ for $\rho_*,\rho^*>0$,
the right-hand side is bounded above by $(e^{-4}-1)\gamma_* a
\downarrow
-\infty$ as $a\uparrow\infty$.
\end{pf}

\subsection{Weak LDP upper bounds}

The weak upper bound for the tagged particle deviations, starting from
local equilibrium measures or
deterministic initial configuration, follows in several steps and is
stated in Step 6. As the same argument works for the current, we also
state its associated weak upper bound in Step 6, below. For the
convenience of the reader, we indicate the modifications needed in
Step 1; the other steps involve similar changes.

\textit{Step 1.}
Consider an interval $[a,b]$ for $0<a<b$; subsequent arguments carry
over straightforwardly to all intervals
$[a,b]\subset\R$
using (\ref{negative_a}) by splitting at the origin if necessary.
Now, divide $[a,b]$
into $m$ equal intervals $A_k=
[c_k,c_{k+1}]$. Then, by the union of events estimate,
\[
\limsup_{N\rightarrow\infty}\frac{1}{N}\log P (X_{N^2T}/N \in
[a,b] ) \leq\max_k\limsup_{N\rightarrow\infty}
\frac{1}{N}\log P ( X_{N^2T}/N \in A_k ).
\]

Then, from (\ref{tagged-current}) and Proposition~\ref{error}, we
have that
\begin{eqnarray*}
&&\limsup_{N\uparrow\infty} \frac{1}{N}\log P (X_{N^2T}/N \in[a,b]
)\\[-1pt]
&&\qquad\leq\limsup_{m\uparrow\infty}
\limsup_{\delta\downarrow
0}\limsup_{n\uparrow\infty}\max_{1\leq k\leq m}\limsup_{N\uparrow
\infty}\frac{1}{N}\\
&& \qquad\quad {}\times \log P \bigl(Y^N_T(G_n) - Y^N_0(G_n)\in
\bigl[Y^N_T\bigl(1_{[0,c_k]}\bigr)-\delta,Y^N_T\bigl(1_{[0,c_{k+1}]}\bigr)+\delta\bigr] \bigr).
\end{eqnarray*}

Since the maps $\mu\mapsto\int G(x) \mu_T\, dx, \int G(x)\mu_0\,dx,
\int_0^c \mu_T\,dx$, for compactly supported $G$ and constants $c$, are
continuous in the Skorohod topology on $D([0,T]; M_1)$,
from Corollary~\ref{KOV}, we
conclude, for fixed $k$, $n$ and $\delta$ that
%
\begin{eqnarray}
\label{upper bound_n_delta}
&&\limsup_{N\rightarrow\infty}\frac{1}{N}\log
p\bigl( Y^N_{T}(G_n) -
Y^N_0(G_n)\,ds
\in\bigl[
Y^N_T\bigl(1_{[0,c_k]}\bigr)-\delta,Y^N_T\bigl(1_{[0,c_{k+1}]}\bigr)
+\delta\bigr]\bigr) \nonumber\hspace*{-15pt}\\
&&\qquad\leq-\inf\biggl\{I_\gamma(\mu); \int
G_n(x) [\mu_{T}(x)-\mu_0(x) ]\,dx  \\
&&\qquad\hphantom{\leq-\inf\biggl\{}
\in\biggl[\int_0^{c_k} \mu_{T}(x)\,dx -\delta,
\int_0^{c_{k+1}}\mu_T(x)\,dx
+\delta\biggr]\biggr\}.\nonumber
\end{eqnarray}

We now indicate the modifications needed for the current in this step.
For $0<a<b$, from (\ref{tagged-current}) and Proposition~\ref{error},
we have
\begin{eqnarray*}
&&\limsup_{N\uparrow\infty} \frac{1}{N}\log P \bigl(J_{-1,0}(N^2T)/N \in
[a,b] \bigr)\\[-1pt]
&&\qquad\leq\limsup_{m\uparrow\infty}
\limsup_{\delta\downarrow
0}\limsup_{n\uparrow\infty}\max_{1\leq k\leq m}\limsup_{N\uparrow
\infty} \frac{1}{N}\log P\bigl (Y^N_T(G_n) - Y^N_0(G_n)\\
&&\qquad \hphantom{\leq\limsup_{m\uparrow\infty}
\limsup_{\delta\downarrow
0}\limsup_{n\uparrow\infty}\max_{1\leq k\leq m}\limsup_{N\uparrow
\infty} \frac{1}{N}\log P\bigl (}\in[c_k-\delta, c_{k+1}+\delta] \bigr).
\end{eqnarray*}
From continuity of the maps $\mu\mapsto\int G \mu_T \,dx$ and
$\mu\mapsto\int G \mu_0\, dx$, and Corollary~\ref{KOV}, we further
bound\vspace*{-1pt}
\begin{eqnarray*}
&&\limsup_{N\rightarrow\infty}\frac{1}{N}\log
P \bigl( Y^N_{T}(G_n) -
Y^N_0(G_n)\,ds \in[c_k-\delta, c_{k+1}+\delta] \bigr) \\[-1pt]
&& \qquad\leq-\inf\biggl\{I_\gamma(\mu); \int
G_n(x) [\mu_{T}(x)-\mu_0(x) ]\,dx \in[c_k-\delta, c_{k+1}+\delta] \biggr\}.
\end{eqnarray*}

\textit{Step 2.} Next, we give a uniform upper bound of the infimum in
(\ref{upper bound_n_delta}). We exhibit
a density $\mu^\c$ satisfying, for each $\delta>0$ and all large $n$,\vspace*{-1pt}
\[
\int G_n(x)[\mu^\c_T(x)-\mu^\c_0(x)]\,dx \in \biggl[\int_0^\c\mu^\c_T(x)\,dx
-\delta, \int_0^\c\mu^\c_T(x)\,dx + \delta\biggr]
\]
and $\sup_{\c\in[a,b]}I_\gamma(\mu^\c)<B_0<\infty$ where $B_0$
is independent of $n$ and $\delta$.

This is accomplished by the constructions in  Section \ref
{finiteness}, namely one takes $\mu^\c= \sigma_t*\gamma+ \lambda
\epsilon(t/T)\psi(x/L)$
with $\lambda, L$ chosen so that
$\lambda L\int_{|\c/L|}^1\psi(x)\,dx
= \int_{u_T}^\c\sigma_T*\gamma(x)\,dx$.
Let $J^\c$ be its current,
and
$H^\c_x$ be the
associated function with respect to (\ref{weak_asym_eq}).

Proposition
\ref{finite_rate_prop}
gives $I_\gamma(\mu^\c)$ is uniformly bounded for $\c\in[a,b]$.
Now compute\vspace*{-1pt}
%
\begin{eqnarray}
\label{G_n_sequence}
&&\int
G_n(x) [\mu^\c_T(x)-\mu^\c_0(x) ]\,dx\nonumber\\
&&\qquad= \int_0^T\int G_n(x)[(1/2)\partial_{xx}\mu^\c-
\partial_x
H^\c_x\mu^\c(1-\mu^\c)]\,dx\,dt \nonumber\\[-8pt]\\[-8pt]
&&\qquad= \int_0^T -(1/2)\partial_x\mu^\c(t,0)
+H^\c_x\mu^\c(1-\mu^\c)(t,0)\,dt\nonumber\\
&&\qquad\quad{} + \frac{1}{n}\int_0^T\int_0^n [ (1/2)\partial
_x\mu^c
-H^\c_x\mu^\c(1-\mu^\c) ]\,dx\,dt.\nonumber
\end{eqnarray}
Since $\int_0^\c\mu^\c_T(x)\,dx = \int_0^T J^\c(0,t)\,dt$ and $J^\c
(0,t) =
-(1/2)\partial_x\mu^\c(t,0) + H^\c_x\mu^\c(1-\mu^\c)(t,0)$,
we have\vspace*{-1pt}
\begin{eqnarray*}
\label{step2bound}
&& \sup_{\c\in[a,b]} \biggl|\int
G_n(x) [\mu^\c_T(x)-\mu^\c_0(x) ]\,dx - \int_0^\c\mu^\c_T(x)\,dx \biggr|
\\
&&\qquad \leq\frac{1}{n} \biggl|\int_0^T\int_0^n (1/2)\partial_x\mu^\c-
H^\c_x\mu^\c(1-\mu^\c) \,dx\,dt \biggr|\\
&&\qquad\leq \sup_{\c\in[a,b]} \biggl|\frac{1}{2n}\int_0^T\bigl(\mu^\c
_t(n)-\mu^\c_t(0)\bigr)\,dt \biggr| +
\frac{1}{\sqrt{n}}\|H^\c_x\mu^\c(1-\mu^\c)\|_{L^2([0,T]\times\R)}.
\end{eqnarray*}
Since $\|H^\c_x\mu^\c(1-\mu^\c)\|^2_{L^2}\leq
2I_0(\mu^\c)$, the right-hand side is $O(n^{-1/2})$ by
Lem\-ma~\ref{finite_rate_prop}.

\textit{Step 3.}
As $I_\gamma$ is a good rate function, by the uniform bounds in step
2, out of minimizers $\nu^{k,n,\delta}$ over
$k=k(m)$, and $n,\delta$ in the infimum in (\ref{upper bound_n_delta}),
by the uniform bound on $I_\gamma(\nu^{k,n,\delta})$,
we can
extract a subsequence, on which the limsup of (\ref{upper
bound_n_delta}) is attained as
$\delta\downarrow0$ and $n,m\uparrow\infty$, and which converges in
$D([0,T]; M_1)$ to a $\bar{\mu}$.

By Proposition~\ref{claudio}, the
subsequence, labeled $\nu^{k,n,\delta}$ itself for simplicity,
may be approximated by $\{\mu^{k,n,\delta}\}$ so
that $\mu^{k,n,\delta}$ is smooth, strictly bounded between $0$
and $1$, $H_x^{k,n,\delta}\in
C^{\infty}([0,T]\times\R)$, Skorohod distance
$d(\mu^{k,n,\delta},\nu^{k,n,\delta})\downarrow0$,
$|I_0(\nu^{k,n,\delta})-I_0(\mu^{k,n,\delta})|\downarrow0$ and when
$\hat{\gamma}\in M_1(\rho_*,\rho^*)$ is piecewise continuous and
$0<\hat{\gamma}(x)<1$ for $x\in\R$,
$|h(\nu^{k,n,\delta}_0; \hat{\gamma})-h(\mu^{k,n,\delta}_0;\hat
{\gamma})|\downarrow
0$.
Also, as
$[a,b]$ is compact, the subsequence can be chosen so that $c_{k+1}$
converges to a $\c\in[a,b]$.

Given $\nu^{k,n,\delta}$ satisfies the restriction in
(\ref{upper bound_n_delta}), we may also arrange
%
\begin{eqnarray}
\label{step3bound}
\int_0^{c_k}\mu^{k,n,\delta}_T(x)\,dx -2\delta&\leq&
\int G_n(x) [\mu^{k,n,\delta}_T(x)-\mu^{k,n,\delta}_0(x)
]\,dx\nonumber\\[-8pt]\\[-8pt]
&\leq& \int_0^{c_{k+1}}\mu^{k,n,\delta}_T(x)\,dx +2\delta.\nonumber
\end{eqnarray}

With these specifications, by lower
semi-continuity, we have that (\ref{upper bound_n_delta}) is
less than, in the case of starting from a local equilibrium measure,
\begin{eqnarray*} &&\lim_{m\uparrow\infty}\lim_{\delta\downarrow
0}\lim_{n\uparrow\infty} \max_k
-I^{LE}_\gamma(\mu^{k,n,\delta})
\leq-I^{LE}_\gamma(\bar{\mu}).
\end{eqnarray*}
When starting from a deterministic configuration, noting
$\nu^{k,n,\delta}_0=\bar{\mu}_0=\gamma$,
(\ref{upper bound_n_delta}) is less than
\[
\lim_{m\uparrow\infty}\lim_{\delta\downarrow
0}\lim_{n\uparrow\infty} \max_k
-I_0(\mu^{k,n,\delta}) \leq-I_0(\bar{\mu})
= -I^{DC}_\gamma(\bar{\mu}).
\]

\textit{Step 4.}
We now show that $\bar{\mu}$ satisfies
%
\begin{equation}
\label{limit_restriction}
\int_0^T \bar{J}(0,t)\,dt = \int_0^\c\bar{\mu}_T(x)\,dx.
\end{equation}

As convergence in $D([0,T]; M_1)$ implies
$\mu^{k,n,\delta}_T\rightarrow\bar{\mu}_T$,
$c_{k+1}-c_k = m^{-1}$ and $0\leq\mu^{k,n,\delta}_T(x)\leq
1$, we have both
\[
\int_0^{c_k}\mu^{k,n,\delta}_T(x)\,dx,
\int_0^{c_{k+1}}\mu^{k,n,\delta}_T(x)\,dx \rightarrow\int_0^\c
\bar{\mu}_T(x)\,dx.
\]
Also, following sequence (\ref{G_n_sequence}),
\begin{eqnarray}
\label{G_neqn}
 &&\int G_n(x) [\mu^{k,n,\delta}_T(x)-\mu^{k,n,\delta}_0(x) ]\,dx\nonumber\\
 && \qquad = \int_0^T J^{k,n,\delta}(0,t)\,dt
\nonumber\\
&&  \qquad  \quad {}+ \frac{1}{n}\int_0^T\int_0^n [(1/2)\partial_x \mu
^{k,n,\delta}_t - H^{k,n,\delta}_x\mu^{k,n,\delta}(1-\mu
^{k,n,\delta})(t,x) ]\,dx\,dt.\nonumber
\end{eqnarray}
Since
$\|H^{k,n,\delta}_x\mu^{k,n,\delta}(1-\mu^{k,n,\delta})\|
^2_{L^2}\leq
2I_0(\mu^{k,n,\delta})$ is uniformly bounded, the last term is
bounded uniformly by
$n^{-1}T + (nT)^{-1/2}\sqrt{2I_0(\mu^{k,n,\delta})}$.
On the other hand,
$\int_0^T J^{k,n,\delta}(0,t)\,dt \rightarrow\int_0^T \bar{J}(0,t)\,dt$
by Proposition~\ref{Lip-J}.

Hence, noting (\ref{step3bound}),
we obtain (\ref{limit_restriction}) immediately.

\textit{Step 5.}
Therefore,
\begin{eqnarray*}
&&\limsup_{m\uparrow\infty}\limsup_{\delta\downarrow
0}\limsup_{n\uparrow\infty}\max_{1\leq k\leq m}\\
&& \quad {}-\inf\biggl\{I_\gamma(\mu);
\int G_n(x) [\mu_{T}(x) -\mu_0(x) ]\,dx\\
&&\hphantom{\quad {}-\inf\biggl\{}\in \biggl[\int_0^{c_k}\mu_{T}(x)\,dx
-\delta,\int_0^{c_{k+1}}\mu_{T}(x)\,dx+\delta\biggr] \biggr\}\\
&& \qquad \leq- I_\gamma(\bar{\mu}) \leq-\min_{\c\in[a,b]}\I(\c).
\end{eqnarray*}

\textit{Step 6.}
The weak LDP upper bound, with respect to the tagged particle, for
compact $K\subset\R$,
%
\begin{equation}
\label{weakldpub}
\limsup_{N\uparrow\infty} \frac{1}{N}P (X_{N^2 t}/N \in K )
\leq-\inf_{a\in K} \I(a),
\end{equation}
is now standard, given that $\I$ is lower
semi-continuous (Lemma~\ref{lowersemicontinuous}).\vadjust{\goodbreak}

Similarly, we have the weak upper bound for the current
%
\begin{equation}
\label{weakldpub_current}
\limsup_{N\uparrow\infty} \frac{1}{N}P \bigl(J_{-1,0}(N^2 t)/N \in K \bigr)
\leq-\inf_{a\in K} \J(a).
\end{equation}

\subsection{LDP lower bound}
\label{lbsection}
As before, we concentrate on the tagged particle deviations, as the
proof for the current is analogous.
For the first step,
the scheme for the weak upper bound is used.
Let $O\subset\R$ be a nonempty open set, and
suppose $a\in O$. We also assume $a>0$ as a similar argument works
for $a\leq0$ by focusing on a subinterval to the left of the origin.
Let $\epsilon>0$ be such that $a-\epsilon>0$ and
$(a-\epsilon,a+\epsilon)\subset O$.

Then, for $\theta>0$,
%
\begin{eqnarray}
\label{lb1}
&& \lim_{N\uparrow\infty}
\frac{1}{N}\log P (X_{N^2 T}/N \in O )\nonumber\\
&& \qquad \geq \lim_{N\uparrow\infty}
\frac{1}{N}P\bigl (X_{N^2 T}/N \in
(a-\epsilon,a+\epsilon) \bigr)\nonumber\\
&& \qquad  \geq\lim_{n\uparrow\infty}\lim_{N\rightarrow\infty}
\frac{1}{N}\log P \Biggl(Y^N_T\bigl(1_{[0,a-\epsilon]}\bigr)<
Y^N_T(G_n)-Y^N_0(G_n)\\
&&\hphantom{\lim_{n\uparrow\infty}\lim_{N\rightarrow\infty}
\frac{1}{N}\log P \Biggl(}\qquad\quad{}+\frac{1}{nN^2}\sum_{x=1}^{nN}J_{x-1,x}(N^2T)<
Y^N_T\bigl(1_{[0,a+\epsilon]}\bigr)\nonumber\\
&&\hspace*{80pt}\hphantom{\lim_{n\uparrow\infty}\lim_{N\rightarrow\infty}
\frac{1}{N}\log P \Biggl(}\qquad\quad \mbox{and }Y^N_T\bigl(1_{[a-\epsilon
,a+\epsilon]}\bigr)>\theta \Biggr).
\nonumber
\end{eqnarray}

From Proposition~\ref{error} and Corollary~\ref{KOV}, the left-hand
side of (\ref{lb1}) is greater
than
%
\begin{eqnarray}\label{lb2}\qquad
&& \lim_{\theta\downarrow0}\lim_{\delta\downarrow
0}\lim_{n\uparrow\infty}\lim_{N\rightarrow
\infty}\frac{1}{N}\log P \bigl(Y^N_T\bigl(1_{[0,a-\epsilon]}\bigr)+\delta<Y^N_T(G_n)-Y^N_0(G_n)
\nonumber\\
&&\hphantom{\lim_{\theta\downarrow0}\lim_{\delta\downarrow
0}\lim_{n\uparrow\infty}\lim_{N\rightarrow
\infty}\frac{1}{N}\log P \bigl(} <Y^N_T\bigl(1_{[0,a+\epsilon]}\bigr)-\delta,
\mbox{ and }Y^N_T\bigl(1_{[a-\epsilon,a+\epsilon]}\bigr)>\theta\bigr)\nonumber\\
&& \qquad\geq \lim_{\theta\downarrow0}\lim_{\delta\downarrow
0}\lim_{n\uparrow\infty} \\
&&\qquad\quad{} -\inf\biggl\{I_\gamma(\mu)\dvtx
\int_0^{a-\epsilon}\mu_T(x)\,dx +\delta
< \int
G_n(x) [\mu_T(x)-\mu_0(x) ]\,dx\nonumber\\
&&\hspace*{44.5pt}\qquad \quad\hphantom{{} -\inf\biggl\{} < \int_0^{a +\epsilon}\mu_T(x)\,dx
-\delta, \mbox{ and }\int_{a-\epsilon}^{a+\epsilon}\mu
_T(x)\,dx>\theta\biggr\}.\nonumber
\end{eqnarray}

Now, for $\alpha>0$, let $\bar{\mu}$ be a density such that
$|I_\gamma(\bar{\mu})-\I(a)|<\alpha$, and
\[
\int_0^T \bar{J}(0,t)\,dt = \int_0^a \bar{\mu}_T(x)\,dx.
\]
By the method used for
(\ref{G_n_sequence}) and (\ref{G_neqn}) in the last section, through
approximations of
$\bar{\mu}$ with smooth $\mu^n$ by Proposition~\ref{claudio}, we
can show
that
%
\begin{equation}
\label{lb_g_n}
\lim_n \int_0^\infty
G_n(x) [\bar{\mu}_T(x) - \bar{\mu}_0(x) ]\,dx = \int_0^T \bar{J}(0,t)\,dt.
\end{equation}

We will need now to approximate $\bar{\mu}$ as follows to ensure a
certain positivity.
Let $\chi= \sigma_s*\gamma+ \lambda\epsilon(t/T)\psi(x/L)$
from  Section~\ref{finiteness} where
$\lambda, L$ are chosen so that
$\int_0^T J^\chi(0,t)\,dt = \int_0^a
\chi_T(x)\,dx$.
Recall $I_\gamma(\chi)<\infty$, and note (\ref{lb_g_n}), with $\chi$
and $J^\chi$ replacing
$\bar{\mu}$ and $J$,
also holds by the explicit construction.
For $0<b<1$, define
$\mu^b = (1-b)\chi+ b\bar{\mu}$. Clearly, $\lim_{b\uparrow1}\mu^b
=\bar{\mu}$ uniformly, and so in $D([0,T]; M_1)$. In fact, $\lim
_{b\uparrow1}
I_\gamma(\mu^b) = I_\gamma(\bar{\mu})$: By lower semi-continuity,
$\liminf
I_\gamma(\mu^b) \geq I_\gamma(\bar{\mu})$ and,
by convexity, $\limsup I_\gamma(\mu^b) \leq I_\gamma(\bar{\mu})$.
Now, for given $\beta>0$, let $b$ be such that $|I_\gamma(\mu^b) -
I_\gamma(\bar{\mu})|<\beta$.

With $\theta>0$, noting
\[
\lim_n \int_0^\infty
G_n(x) [\mu^b_T(x) - \mu^b_0(x) ]\,dx = \int_0^a \mu^b_T(x)\,dx,
\]
we have for $n\geq N(\theta, \bar\mu, \chi)$ that
\begin{eqnarray*}
&&\int_0^{a-\epsilon} \mu^b_T(x)\,dx +b\int_{a-\epsilon}^a \bar{\mu}_T(x)\,dx
+(1-b)\int_{a-\epsilon}^a\chi_T(x)\,dx - \theta\\
&&\qquad\leq\int G_n(x) [\mu^b_T(x)-\mu^b_0(x) ]\,dx\\
&&\qquad\leq\int_0^{a+\epsilon} \mu^b_T(x)\,dx
-b\int_a^{a+\epsilon}\bar{\mu}_T(x)\,dx -
(1-b)\int_a^{a+\epsilon}\chi_T(x)\,dx
+ \theta.
\end{eqnarray*}
By the construction of $\chi$,
$\int_{a-\epsilon}^a\chi_T(x)\,dx, \int_a^{a+\epsilon}\chi_T(x)\,dx
\geq
c\epsilon$ for a constant $c>0$.
Hence,
we can choose $\theta=\theta(\epsilon, b, \chi)$ so that for all
small $\delta$,
\[
(1-b)\int_a^{a+\epsilon}\chi_T(x)\,dx-\theta, (1-b)\int_{a-\epsilon
}^a\chi_T(x)\,dx -\theta
> \delta.
\]

Therefore, as $\bar{\mu}$ is nonnegative, $\mu^b$
satisfies the restriction in the infimum in
(\ref{lb2}). In particular,
we have
\[
\lim_{N\uparrow\infty} \frac{1}{N}P
(X_{N^2 t}/N \in O )
\geq-I_\gamma(\mu^b)\geq-I_\gamma(\bar{\mu}) -\beta\geq-\I(a)
- \alpha-\beta.
\]
Hence,
%
\begin{equation}
\label{lowerboundldp}
\lim_{N\uparrow\infty} \frac{1}{N}P
(X_{N^2 t}/N \in O ) \geq-\inf_{a\in O} \I(a).
\end{equation}

Analogously, we have weak lower bound large deviations for the current
%
\begin{equation}
\label{lowerboundldp_current}
\lim_{N\uparrow\infty} \frac{1}{N}P
\bigl(J_{-1,0}(N^2 t)/N \in O \bigr) \geq-\inf_{a\in O} \J(a).
\end{equation}

\begin{pf*}{Proof of Theorem~\ref{mainthm1}}
First, the functions $\J$ and $\I$ are finite-valued rate functions
which vanish exactly at $v_T$ and $u_T$, respectively, by Corollary~\ref{I_ratefunction}.

Next, a ``weak'' LDP is found from (\ref{weakldpub_current}) and (\ref
{lowerboundldp_current}) with respect to rate function~$\J$, and (\ref
{weakldpub}) and
(\ref{lowerboundldp}) with respect to $\I$.
Standard arguments, given exponential
tightness (Lemma~\ref{exponentialtightnessLemma}), extend the ``weak'' LDP
to the full large deviation principle.

Finally, given the LDP and exponential tightness, it follows that (1)~$\J$ and $\I$ are good rate functions by Lemma
1.2.18~\cite{DZ}, and also that (2) $\lim_{|a|\uparrow\infty} \J
(a) = \lim_{|a|\uparrow\infty}\I(a)=\infty$.
\end{pf*}

\section{Asymptotic evaluations}\label{asymptotic_eval_section}

We prove Theorems~\ref{mainthm2},~\ref{mainthm3} and~\ref{mainthm4}
in succeeding subsections.

\subsection{\texorpdfstring{Proof of Theorem \protect\ref{mainthm2}}{Proof of Theorem 1.6}} We first prove the
upper bounds which are implied by the following lemma, and then the
lower bounds.

\begin{lemma}
\label{upper bound_asymp}
Starting from (DIC) or (LEM) initial conditions, there are constants
$c_0,c_1,c_2,c_3$ depending only on $\gamma$, such that when,
respectively, $|a-v_T|/\sqrt{T}\geq c_0$ and $|a-u_T|/\sqrt{T} \geq
c_0$, we have, in turn,
\[
\J(a) \leq\frac{c_1|a-v_T|^3}{T} \quad\mbox{and}\quad\I(a) \leq
\frac{c_1|a-u_T|^3}{T}.
\]
Also, when, respectively, $|a-v_T|/\sqrt{T}\leq c_2$ and
$|a-u_T|/\sqrt{T}\leq c_2$, we have, correspondingly,
\begin{eqnarray*}
\J(a) \leq \frac{c_3(a-v_T)^2}{\sqrt{T}} \quad\mbox{and}\quad\I
(a) \leq \frac{c_3(a-u_T)^2}{\sqrt{T}}.
\end{eqnarray*}
\end{lemma}

\begin{pf} We prove the estimates for the current rate function, and
deduce corresponding bounds for the tagged particle rate function.
Let also $a>v_T$ as the argument for $a<v_T$ is analogous.
For the reader's convenience, we recall estimate (\ref{I_0_bound}) and write
\[
\J(a) \leq\frac{4\epsilon^*}{\gamma_*(1-\gamma^*)} \biggl[\frac
{\lambda^2 T}{4
|L|}\int_{-1}^1\psi'(x)^2\,dx +
\frac{\lambda^2|L|^3}{T}\int_{-1}^1\Psi(x)^2\,dx \biggr].
\]
Recall also the restriction equation (\ref{I_restriction}) when $\c= a$,
\[
\lambda= \frac{a-v_T}{L\int_0^1\psi\, dx},
\]
subject to $0<\lambda\leq\min\{\gamma_*,1-\gamma^*\}/2$. The
requirement on $\lambda$ holds
when
\[
L \geq|a-v_T|\Big/\biggl[\frac{1}{2}\min\{\gamma_*,(1-\gamma^*)\}\int
_{0}^1\psi dx\biggr] := \kappa_0|a-v_T|.
\]
Now take $L$ in the form $L = \kappa\sqrt{T}$.
Substituting into
the bound for $\J(a)$,
we obtain
\[
\J(a) \leq\frac{(a-v_T)^2}{\sqrt{T}}\frac{4\epsilon^*\kappa
}{\gamma_*(1-\gamma^*)} \biggl[\frac{1}{4\kappa^4}\int_{-1}^1\psi'(x)^2\,dx+
\int_{-1}^1\Psi(x)^2\,dx \biggr].
\]

Hence, when $a$ is large, say
$\kappa = |a-v_T|\kappa_0/\sqrt{T} \geq 1$,
we have
$\J(a) \leq c(\gamma)|a-v_T|^3/T$.
Correspondingly, when $a$ is such that $|a-v_T|\kappa_0/\sqrt{T}\leq
1$, we choose $\kappa= 1$ to get
$\J(a) \leq c(\gamma)|a-v_T|^2/\sqrt{T}$.

The bounds on the tagged particle rate function $\I$ follow from the
current rate function bounds. First, by (\ref{initial_ineq}), $\I
^{LE}_\gamma(a) \leq\I_\gamma^{DC}(a)$. Also, by (\ref
{DIC_relation}), with $\gamma'(x) = \gamma(x+a)$,
$\I^{DC}_\gamma(a) = \J_{\gamma'}^{DC}(\int_0^a \gamma \,dx)$. For
fixed $a$, let now $v_T(\gamma')$ be the LLN integrated current
through the origin starting from $\gamma'$. Then
\begin{eqnarray*}
\int_0^a \gamma dx - v_T(\gamma') & = & \int_0^{a}\gamma \,dx - \int
_0^\infty\sigma_T* \gamma' - \gamma'\,dx \\[-2pt]
&=& \int_0^a \sigma_T *\gamma \,dx - \int_0^\infty\sigma_T * \gamma
- \gamma \,dx= \int_{u_T}^a \sigma_T*\gamma \,dx.
\end{eqnarray*}
Hence, $\gamma_*|a-u_T| \leq| \int_0^a \gamma \,dx - v_T(\gamma
')| \leq\gamma^*|a-u_T|$. Since $\gamma_*, \gamma^*$ are
uniform lower and upper bounds on $\sigma_{T/10}*\gamma$ (and hence
on $\sigma_{T/10}*\gamma'$), the desired estimates on $\I(a)$ are
derived from the bounds on $\J_{\gamma'}(\int_0^a \gamma \,dx)$.
\end{pf}

The lower bounds in Theorem~\ref{mainthm2} are implied by the
following two estimates.

\begin{lemma}\label{lowerbound1} Starting from a (DIC) condition with
profile $\gamma$, there are constants $c_0=c_0(\gamma
,T),c_1=c_1(\gamma)$ such that for $|a| \geq c_0$, we have
\[
\J(a), \I(a) \geq\frac{c_1|a|^3}{T}.
\]
\end{lemma}

\begin{lemma}\label{lowerbound2} Starting from a (DIC) condition with
profile $\gamma$, there is a constant $c_1(\gamma)$, such that for
$a\in\R$, we have
\[
\J(a) \geq\frac{c_1(a-v_T)^2}{\sqrt{T}} \quad\mbox{and}\quad\I
(a) \geq\frac{c_1(a-u_T)^2}{\sqrt{T}}.
\]
\end{lemma}

\begin{pf*}{Proof of Lemma~\ref{lowerbound1}} We concentrate first on the
current calculation. Suppose $a>0$, as the argument for $a<0$ is
similar. Let $\hat\gamma\in M_1(\rho_*,\rho^*)$ be a smooth
density, strictly bounded away from $0$ and $1$, such that $h(\gamma;
\hat\gamma)<\infty$. For
$\epsilon>0$, by Proposition~\ref{claudio}, let $\mu$ be a smooth
density such that $\mu_0= \sigma_\alpha*\gamma$, $|h(\mu_0; \hat
\gamma)-h(\gamma; \hat\gamma)|<\epsilon$ and
$|I_0(\mu) - \J(a)| \leq \epsilon$.
Noting Proposition~\ref{Lip-J}, we can, in addition, impose on the
approximating density that
$
|\int_0^TJ(0,t)\,dt - a| \leq \epsilon$.

Now, noting (\ref{current-massformula}) in Proposition \ref
{energyLemma}, we have the Lipschitz bound,
\[
\biggl|\int_0^TJ(x,t)\,dt - \int_0^TJ(0,t)\,dt \biggr| = \biggl|\int_0^x
\mu_T(z) - \mu_0(z)\,dz \biggr| \leq|x|.\vadjust{\goodbreak}
\]
Then,
for $0\leq x\leq a - \epsilon$, we have
$\int_0^T J(x,t)\,dt \geq \int_0^TJ(0,t)\,dt -x \geq a -\epsilon-x$
so that
\begin{eqnarray*}
(a - \epsilon)^3/3 &=& \int_0^{a - \epsilon}[a - \epsilon- x]^2\,dx \\
&\leq& \int_0^{a-\epsilon} \biggl[\int_0^T J(x,t)\,dt \biggr]^2\,dx \leq T\int\int
_0^T J^2\,dt\,dx.
\end{eqnarray*}
Hence, as $\mu(1-\mu), \hat\gamma(1-\hat\gamma)\leq1/4$, from
the formula for $I_0(\mu)$ in Proposition~\ref{energyLemma} and
simple computations,
%
\begin{eqnarray}\label{large_a}\qquad
\J(a)\geq I_0(\mu)-\epsilon& \geq& \int\int_0^T J^2\,dt\,dx - \frac
{1}{2}h(\mu_0;\hat\gamma)
-T\int\frac{(\partial_x\hat\gamma)^2}{\hat\gamma^2(1-\hat\gamma
)^2}\,dx - \epsilon\nonumber\\[-8pt]\\[-8pt]
&\geq&\frac{(a-\epsilon)^3}{3T} - \frac{1}{2}h(\gamma;\hat\gamma)
-T\int\frac{(\partial_x\hat\gamma)^2}{\hat\gamma^2(1-\hat\gamma
)^2}\,dx- \frac{3}{2}\epsilon.\nonumber
\end{eqnarray}

For the tagged particle rate function, from (\ref{DIC_relation}), we
have $\I^{DC}_\gamma\hspace*{-0.6pt}(a)=\break
\J^{DC}_{\gamma'}\hspace*{-0.6pt}(\int_0^a\hspace*{-0.6pt}\gamma dx)$
where $\gamma'(x) = \gamma(x+a)$. Since $\gamma(x) \geq\min\{\rho
_*,\rho^*\}$ for all large $|x|$, $|\int_0^a \gamma dx| \geq c(\gamma
)|a|$ for all large $|a|$ where $c(\gamma)>0$. Also, as $\gamma',
\hat\gamma\in M_1(\rho_*,\break\rho^*)$, by calculation $h(\gamma'; \hat
\gamma) = O(|a|)$. Hence, plugging into (\ref{large_a}), we obtain
the desired estimate on $\I(a)$.
\end{pf*}

\begin{pf*}{Proof of Lemma~\ref{lowerbound2}} We focus first on the current
rate function computation.
By Proposition~\ref{claudio} and~\ref{Lip-J}, let $\mu$ be a smooth
density with properties (i)--(viii) such that $\mu_0 = \sigma_\alpha
*\gamma$, $|\J(a) - I_0(\mu)|<\epsilon$ and $|\int_0^T J(t,0)\,dt -
a|<\epsilon$.
Let $v_T(\alpha)$ be the LLN speed starting from profile $\sigma
_\alpha*\gamma$, and note $\lim_{\alpha\downarrow0}|v_T -
v_T(\alpha)|=0$.

Consider solutions of $\partial_t\rho= (1/2)\rho_{xx}$ and $\partial
_t\mu= (1/2)\partial_{xx}\mu-\partial_x(H_x\mu(1-\mu))$, both
with initial value $\sigma_\alpha*\gamma$. The difference $U = \rho
- \mu$ satisfies equation $\partial_tU = (1/2)\partial_{xx}U -
\partial_x(H_x\mu(1-\mu))$ with $U(0,x) \equiv0$. Integrating once
in the space variable, noting properties of $\mu$, $S(t,x) = \int
_{-\infty}^x U(t,y)\,dy$ satisfies
$\partial_t S = (1/2)\partial_{xx}S - H_x\mu(1-\mu)$. Hence,
we have
\begin{eqnarray*}
S(t,x) &=& \sigma_t*S(0,x)\\
&& {}+ \int_0^t\int\frac{1}{\sqrt{2\pi(t-s)}}e^{\fracc
{-(x-y)^2}{2(t-s)}}[-H_x\mu(1-\mu)](s,y)\,dy\,ds\\
&=& -\int_0^t\int\frac{1}{\sqrt{2\pi t}} e^{\fracc{-(x-y)^2}{2t}}
H_x\mu(1-\mu)(t-s,y)\,dy\,ds.
\end{eqnarray*}

Now, the difference in integrated macroscopic currents across $x$ up to
time $t$ with respect to $\rho$ and $\mu$ is
$-S(t,x)$; cf. above (\ref{converges_eqn}). Therefore, by the Schwarz
inequality and $0\leq\mu\leq1$, when $x=0$, we have for small
$\alpha$ that
\begin{eqnarray*}
&&\bigl(v_T - a +O(\epsilon)\bigr)^2 \\
&& \qquad\leq\biggl[\int_0^T \biggl(\int\sigma^2_t(y)\,dy
\biggr)^{1/2} \biggl(\int H_x^2\mu(1-\mu)(y,t)\,dy\biggr)^{1/2}\,dt\biggr]^2.
\end{eqnarray*}
As $\|\sigma^2_t\|^2_{L^2(\R)} \leq C t^{1/2}$, a further bound of
the right-hand side is $2C\sqrt{T}I_0(\mu) \leq2C\sqrt{T}(\J(a) +
O(\epsilon))$ for some universal constant $C$.

We now use relations (\ref{DIC_relation}) to analyze the tagged
particle rate function. Indeed, let $u_T(\alpha)$ be the corresponding
tagged particle LLN speed starting from profile $\rho_0 = \sigma
_\alpha*\gamma$, and note $\lim_{\alpha\downarrow0}u_T(\alpha
)=u_T$. As before, by Proposition~\ref{claudio}, let $\mu$ be a
smooth density such that $\mu_0=\rho_0$, $|\I(a) - I_0(\mu
)|<\epsilon$ and by Proposition~\ref{Lip-J}, $|\int_0^T J(a,t)\,dt -
\int_0^a \rho_0(x)\,dx|<\epsilon$.

Note, with respect to density $\rho$, the current across $a$
equals $\int_a^{\infty}\sigma_T*\rho_0 - \rho_0 \,dx$, and the current
across the origin equals $\int_0^\infty\sigma_T*\rho_0 - \rho_0 \,dx
= \int_0^{u_T(\alpha)}\sigma_T*\rho_0 \,dx$. Then, for small $\alpha$,
the square of the difference in integrated currents with respect to
$\rho$ and $\mu$ across $a$ equals $(
\int_{u_T}^a\sigma_T*\gamma(x)\,dx+O(\epsilon))^2 \geq
\gamma_*^2(a-u_T + O(\epsilon))^2$ where $\sigma_T*\gamma\geq
\gamma_*>0$. But, on the other hand, as before, $\|\sigma_t(\cdot
-a)\|^2_{L^2(\R)}\leq Ct^{1/2}$, and the square current difference
is still bounded by $2C\sqrt{T}I_0(\mu) \leq2C\sqrt{T}(\I(a) +
O(\epsilon))$. This finishes the proof.
\end{pf*}

\subsection{\texorpdfstring{Proof of Theorem \protect\ref{mainthm3}}{Proof of Theorem 1.7}}
Starting from a (DIC) state, since $\gamma(x)\equiv\rho$, noting
(\ref{DIC_relation}), we observe that $\gamma'(x) \equiv\rho$,
$\int_0^a \gamma \,dx = a\rho$, and $\I(a) = \J(a\rho)$. Hence, we
need only give the argument for the current rate function $\J$, as the
estimate for the tagged particle rate function $\I$ follows directly.

We now make some useful reductions. Recall, when starting under a
deterministic configuration with profile
$\gamma(x)\equiv\rho$, in order for $I^{\mathrm{DC}}_\gamma
(\nu)<\infty$,
$\nu$
must satisfy $\nu_0(x) \equiv\rho$ and $I^{\mathrm{DC}}_\gamma
(\nu) = I_0(\nu)<\infty$.
By Proposition~\ref{claudio} and Proposition~\ref{Lip-J}, for each
$\epsilon>0$, we can find a smooth density $\mu$,
such that $\mu_0(x)\equiv\rho$ and
\[
\J(a) \geq I_0(\mu) - \epsilon\frac{a^2}{\sqrt{T}}.
\]
In addition, we may also impose that
\[
\biggl|\int_0^T J(0,t)\,dt - a \biggr|< a\epsilon.
\]

For such a density $\mu$,
by
Proposition~\ref{energyLemma} [applied with $\gamma(x)\equiv\rho$],
%
\begin{eqnarray}
\label{rate-relation1}
I_0(\mu) &=& \frac{1}{8}\int\int_0^T \frac{(\partial_x\mu)^2}{\mu
(1-\mu)}\,dt\,dx + \frac{1}{2}h(\mu_T;\rho) \nonumber
\\[-8pt]
\\[-8pt]&&{}+ \frac{1}{2}\int\int
_0^T\frac{J^2}{\mu(1-\mu)}\,dt\,dx.
\nonumber
\end{eqnarray}

Consider now a sequence $\{\mu^a\}$ of such
$\epsilon a^2/\sqrt{T}$-minimizers of
$\J(a)$ as \mbox{$|a|\downarrow0$}. The upper and lower bounds
in Lemmas~\ref{upper bound_asymp} and~\ref{lowerbound2}, as $v_T=0$,
gives $I_0(\mu^a)=O(a^2/\sqrt{T})$.
Then, by Lemma
\ref{unif_limit}, we have
$\mu^a\rightarrow\rho$ in $L^2([0,T]\times\R)$, and in fact
\[
\sup_{0\leq t\leq T}\int\bigl(\mu^a(t,x)-\rho\bigr)^2\,dx = O(a^2).
\]

We now deduce that there are functions $r(t,x)$ and $j(t,x)$ on
$[0,T]\times\R$ such that $r(0,x)\equiv0$, $\partial_t r + \partial
_x j = 0$ weakly in $L^2([0,T]\times\R)$, $\int_0^T j(t,0)\,dt = 1$, and
%
\begin{eqnarray}
\label{weak_relation}\qquad
&&\rho(1-\rho)\times\liminf_{a\downarrow0} \J(a)/a^2 \nonumber\\
[-8pt]\\[-8pt]
&&\qquad\geq\frac{1}{8}\int_0^T \int(\partial_x r)^2 \,dx\,dt + \frac
{1}{4}\int|r(T,x)|^2\,dx + \frac{1}{2}\int_0^T\int
j^2(t,x)\,dx\,dt.\nonumber
\end{eqnarray}

Consider a function $\lambda^a(t,x) = \psi(\mu^a(t,x))$ where $\psi
'(x) = \min\{(x(1-x))^{-1/2},\break  M\}$ for some $M\geq2(\rho(1-\rho
))^{-1/2}$. Then, $\partial_x\lambda^a = \psi'(\mu^a)\partial_x\mu
^a \leq(\mu^a(1-\mu^a))^{-1/2}\partial_x\mu^a$, and so
\[
\int_0^T\int(\partial_x \lambda^a)^2 \,dx\,dt \leq\int_0^T\int\frac
{(\partial_x\mu^a)^2}{\mu^a(1-\mu^a)}\,dx\,dt.
\]

At this point, let us take weak $L^2([0,T]\times\R)$ limits of
$a^{-1}\partial_x \lambda^a$, $a^{-1}(\mu^a - \rho)$ and
$a^{-1}J^a$, and label them as $u$, $r$ and $j$, respectively. Also,
take a weak $L^2(\R)$ limit of $a^{-1}(\mu^a(T,x) -\rho)$ and call
it $q$. Using suitable truncations, and Fatou's Lemma, given $\mu^a
\rightarrow\rho$ strongly, we have
\begin{eqnarray*}
\int_0^T \int u^2\,dx\,dt & \leq& \liminf\frac{1}{a^2}\int_0^T\int
(\partial_x \lambda^a)^2 \,dx\,dt,\\
\frac{1}{2\rho(1-\rho)}\int|q(x)|^2\,dx & \leq& \liminf\frac
{1}{a^2}\int h_d(\mu^a_T(x);\rho)\,dx,\\
\frac{1}{\rho(1-\rho)}\int_0^T\int j^2\,dx\,dt & \leq& \liminf\frac
{1}{a^2}\int_0^T\int\frac{(J^a)^2}{\mu^a(1-\mu^a)}\,dx\,dt.
\end{eqnarray*}

We may identify (a) $\partial_x r = \sqrt{\rho(1-\rho)}u$, (b)
$r(T,x) = q(x)$, and (c) $\partial_t r + \partial_x j = 0$ weakly in
$L^2([0,T]\times\R)$. The last two (b), (c) follow from weak limits and
properties of $\mu^a$. However, (a) also holds given the weak limits
since $\partial_x \mu^a = \psi'(\mu^a)^{-1}\partial_x\lambda^a$
and $\psi'(\mu^a)^{-1} \rightarrow\sqrt{\rho(1-\rho)}$ strongly
in $L^2$.

Now,
define
\[
K(t,x) = \int_0^t j(s,x)\,ds.
\]
Then, the right-hand side of (\ref{weak_relation}) becomes
\begin{eqnarray*}
\mathcal K & = & \frac{1}{4}\int|\partial_xK(T,x)|^2\,dx + \frac
{1}{2}\int_0^T\int|\partial_tK(t,x)|^2\,dx\,dt\\
&&{}+ \frac{1}{8}\int_0^T\int|\partial_{xx}K(t,x)|^2\,dx\,dt.
\end{eqnarray*}
By scaling, $M(t,x) = K(tT,x\sqrt{T})$, we obtain
\[
\liminf_{|a|\downarrow0} \frac{\sqrt{T}}{a^2}\J(a) \geq[\rho
(1-\rho)]^{-1} \inf\mathcal M,
\]
where the infimum is over $M\in C^{1,2}([0,1]\times\R)$, such that
$M(0,x)\equiv0$ and $M(1,0)=1$, and
\begin{eqnarray*}
\mathcal M & = & \frac{1}{4}\int|M_x(1,x)|^2\,dx + \frac{1}{2}\int
_0^1\int|M_t(t,x)|^2\,dx\,dt \\
&&{}+ \frac{1}{8}\int_0^1\int|M_{xx}(t,x)|^2\,dx\,dt.
\end{eqnarray*}

On the other hand, the upper bound
%
\begin{equation}
\label{var_upper bound}
\limsup_{|a|\downarrow0} \frac{\sqrt{T}}{a^2}\mathbb J(a) \leq[\rho
(1-\rho)]^{-1}\inf\mathcal M
\end{equation}
also follows by a similar strategy:
In Proposition~\ref{mathcal_M} below, we evaluate $\inf\mathcal M$
and find a minimizer. One can find a smooth $\epsilon$, approximating
$M$ with bounded derivatives, and trace back to obtain the
corresponding density $\mu^a$ satisfying $a^{-1}(\mu^a-\rho) =
\partial_x K$, $a^{-1}J^a = \partial_t K$, $a^{-1}\partial_x \mu^a
= \partial_{xx}K$ with $\int_0^TJ^a(0,t)\,dt = a$ and $\mu
^a_0(x)\equiv\rho$.
Given $\|\partial_x M\|_{L^\infty([0,T]\times\R)}<\infty$, we have
$\|\mu^a - \rho\|_{L^\infty([0,T]\times\R)} \leq|a|\|\partial_x
K\|_{L^\infty} = (|a|/\sqrt{T})\|\partial_x
M\|_{L^\infty} = O(|a|)$.
The argument to derive (\ref{var_upper bound}) now follows from
standard approximations with respect to (\ref{rate-relation1}).

Hence, the proof of Theorem~\ref{mainthm3} will follow
from evaluations $\inf_M \mathcal M = \sqrt{\pi}/2$,
$\sigma_{X, dyn}^2= (1-\rho)/(\rho\sqrt{\pi})$ and $\sigma_{J,
dyn}^2 = \rho(1-\rho)/\sqrt{\pi}$ in Propositions
\ref{mathcal_M} and~\ref{dynamical_part} below.

\begin{proposition} \label{mathcal_M}
We have
\[
\inf_M \mathcal M = \frac{\sqrt{\pi}}{2},
\]
where the infimum is over $M\in C^{1,2}([0,1]\times\R)$ such that
$M(0,x)\equiv0$ and $M(1,0)= 1$.
\end{proposition}

\begin{pf} The argument is in three steps.
(A) We first minimize
%
\begin{equation}\label{double_star}
\int_0^1\int_{-\infty}^\infty\frac{1}{2} |M_t(t,x)|^2 + \frac{1}{8}
|M_{xx}(t,x)|^2\,dx\,dt
\end{equation}
when $M(0,x)\equiv0$ and $M(1,x)$ is a given compactly supported
$C^4(\R)$ function.
The Euler equation is
%
\begin{equation}
\label{one_star}M_{tt} = \frac{1}{4}M_{xxxx}
\end{equation}
with the boundary conditions at $t=0,1$.

One can verify the solution of \eqref{one_star}, which is smooth and
classical,
in terms of Fourier transform with respect to the $x$ variable but not
transforming the $t$ variable, is given by
%
\begin{equation}
\label{no_star}
\hat{M}(t,y) = \hat{M}(1,y)\frac{e^{\fraca{ty^2}{2}} -e^{\fraca
{-ty^2}{2}}}{e^{\fraca{y^2}{2}} - e^{\fraca{-y^2}{2}}},
\end{equation}
where
\[
\hat{M}(1,y) = \frac{1}{\sqrt{2\pi}}\int e^{iyx}M(1,x)\,dx.
\]
The corresponding value of \eqref{double_star}, through Plancherel's
formula, is expressed as
\[
\int_{-\infty}^\infty|\hat{M}(1,y)|^2k(y)\,dy,
\]
where
\begin{eqnarray*}
k(y)&=& \int_0^1 \frac{y^4}{8} \frac{ [e^{\fraca{ty^2}{2}}+e^{\fraca
{-ty^2}{2}} ]^2 + [e^{\fraca{ty^2}{2}}-e^{\fraca{-ty^2}{2}} ]^2}
{ [e^{\fraca{y^2}{2}} -e^{\fraca{-y^2}{2}} ]^2}\,dt\\
&=&\int_0^1\frac{y^4}{4}\frac{e^{ty^2} +e^{-ty^2}}{ [e^{y^2/2}
-e^{-y^2/2} ]^2}\,dt\\
&=&\frac{y^2}{4}\frac{e^{y^2}-e^{-y^2}}{ [e^{y^2/2} -e^{-y^2/2} ]^2}=
\frac{y^2}{4}\frac{e^{y^2/2} +e^{-y^2/2}}{e^{y^2/2}-e^{-y^2/2}}.
\end{eqnarray*}
Given that the integrand in \eqref{double_star} is a strict convex function
of $M_t$ and $M_{xx}$,
solution \eqref{no_star} is the unique minimizer of
\eqref{double_star} (by say straightforward modifications
of the
proof of~\cite{DAC}, Theorem 2.1).

(B) Now, we consider
the term
\[
\frac{1}{4}\int_{-\infty}^\infty|M_x(1,x)|^2\,dx = \frac{1}{4}\int
_{-\infty}^\infty y^2|\hat{M}(1,y)|^2\,dy
\]
and
minimize, over $M\in L^2([0,1]\times\R)$,
\begin{eqnarray*}
&&\frac{1}{4}\int_{-\infty}^\infty|\hat{M}(1,y)|^2 y^2 \biggl[ 1+ \frac
{e^{y^2/2} +e^{-y^2/2}}{e^{y^2/2}-e^{-y^2/2}} \biggr]\,dy\\
&&\qquad= \frac{1}{2}\int_{-\infty}^\infty|\hat{M}(1,y)|^2 y^2
\frac{e^{y^2/2}}{e^{y^2/2}-e^{-y^2/2}}\,dy
\end{eqnarray*}
subject to
\[
\frac{1}{\sqrt{2\pi}}\int_{-\infty}^\infty\hat{M}(1,y)\,dy = 1.
\]

Recall that the minimizer of
\[
\int|g(y)|^2 K(y)\,dy
\]
when $\int g(y)\,dy = a$ is given by $g(y) = cK(y)^{-1}$ and $c=a[\int
K(y)^{-1}\,dy]^{-1}$, with minimum value $a^2[\int K(y)^{-1}\,dy]^{-1}$.
Hence, with $a=\sqrt{2\pi}$ and
\[
K(y) = \frac{y^2}{2}\frac{e^{y^2/2}}{e^{y^2/2}-e^{-y^2/2}},
\]
we identify $M(1,y)$ through its transform $\hat M(1,y) = cK(y)^{-1}$.
Denote $\widetilde M$ as the function in \eqref{no_star} with this
choice of $M(1,y)$.

(C) Let now $M^*$ be a compactly supported $C^{2,4}([0,1]\times\R)$
function such that $M^*(0,x)\equiv0$ and $M^*(1,0)=1$ whose $\mathcal
M$-value approximates $\inf_M \mathcal M$. From steps (A) and (B), we
obtain a lower bound of the infimum value which is actually achieved by
the smooth $C^{2,4}([0,1]\times\R)$ function $\widetilde M$.
Therefore, $\widetilde M$ is a minimizer.

Finally, given
\begin{eqnarray*}
\int K(y)^{-1}\,dy & = & 2\int\frac{1-e^{-y^2}}{y^2}\,dy\\
& = & 2\int_0^1\int e^{-ty^2}\,dy\,dt = 2 \int_0^1\sqrt{\frac{\pi
}{t}}\,dt = 4\sqrt{\pi},
\end{eqnarray*}
we obtain the infimum,
$\inf_M \mathcal M = 2\pi/(4\sqrt{\pi}) = \sqrt{\pi}/2$, as
desired.
\end{pf}

\begin{proposition}
\label{dynamical_part}
Starting under initial distribution $\nu_\rho$, the dynamical parts
of the limiting variances of $T^{-1/4}J_{-1,0}(T)$ and $T^{-1/4}x(T)$ under
$\nu_\rho$
are
\[
\sigma^2_{J, dyn} := \lim_{T\rightarrow\infty}
\frac{1}{\sqrt{T}}E_{\nu_\rho} \bigl[\bigl (J_{-1,0}(T) -
E_\eta[J_{-1,0}(T)] \bigr)^2 \bigr] = \frac{\rho(1-\rho)}{\sqrt{\pi}}
\]
and
\[
\sigma^2_{X, dyn} := \lim_{T\rightarrow\infty}
\frac{1}{\sqrt{T}}E_{\nu_\rho} \bigl[\bigl (x(T) -
E_\eta[x(T)] \bigr)^2 \bigr] = \frac{1-\rho}{\rho}\frac{1}{\sqrt{\pi}}.\vadjust{\goodbreak}
\]
\end{proposition}

\begin{pf}
First, we note the limit distribution and
variance
of both\break $T^{-1/4}x(T)$ and $\rho^{-1}T^{-1/4}J_{-1,0}(T)$
are the same, namely $N(0,\sigma^2)$ with\break
$\sigma^2=\sqrt{2/\pi}(1-\rho)/\rho$; cf.~\cite{Arratia}.
Moreover, $\lim_{T\uparrow\infty}E_{\nu_\rho}[(T^{-1/4}X(T) -\break \rho
^{-1}T^{-1/4}J_{-1,0}(T))^2]=0$, since $(T^{-1/4}X(T) - \rho
^{-1}T^{-1/4}J_{-1,0}(T))^2$ vanishes in probability, and also is
uniformly integrable; cf.~\cite{DeMF}, equation (28), or \cite
{Arratia}, page 368, and~\cite{Peligrad-S}, Proposition 4.2 and proof of
Lemma 3.2.

Then we need only show
\[
\lim_{T\rightarrow\infty} \frac{1}{\sqrt{T}}E_{\nu_\rho} [
(E_\eta[J_{-1,0}(T)] )^2 ] = \rho(1-\rho)\frac{\sqrt{2}-1}{\sqrt
{\pi}},
\]
which, given the form of the limiting variance of the scaled current, and
\[
E_{\nu_\rho}[(J_{-1,0}(T))^2] =E_{\nu_\rho}\bigl[\bigl(J_{-1,0}(T)-E_\eta
[J_{-1,0}(T)]\bigr)^2\bigr] +
E_{\nu_\rho}[(E_\eta[J_{-1,0}(T)])^2]
\]
implies the desired results.

Now, the current $J_{-1,0}$ has martingale decomposition (cf. Section 2
\cite{Peligrad-S}),
\[
J_{-1,0}(t) = M(t) + \frac{1}{2}\int_0^t \eta_s(-1)-\eta_s(0)\,ds.
\]
Also, for $x\in\Z$, from ``duality'' (cf. Liggett
\cite{Liggett}, Section VIII.1, page 363),
\[
E_\eta[ \eta_t(x) ] = \sum_i p(t,i-x) \eta(i),
\]
where $p(t,j)=P(S_t = j)$ is the probability a continuous time random
walk, starting from the origin, travels to $j$ in time $t$.
Then,
\begin{eqnarray*}
E_\eta[J_{-1,0}(T)] & = & \frac{1}{2}\int_0^T E_\eta[\eta_t(-1)] -
E_\eta[\eta_t(0)]\,dt\\[-2pt]
& = & \frac{1}{2}\sum_i \eta(i)\int_0^T p(t,i+1)-p(t,i)\,dt\\[-2pt]
&=&\frac{1}{2}\sum_i \bigl(\eta(i)-\rho\bigr)\int_0^T p(t,i+1)-p(t,i)\,dt.
\end{eqnarray*}
Therefore, from independence of coordinates $\{\eta(i)\}$,
\[
Q_0(T):= E_{\nu_\rho} [ (E_\eta[J_{-1,0}(T)] )^2 ] = \rho(1-\rho)
\sum_i \biggl| \frac{1}{2}\int_0^Tp(t,i+1)-p(t,i)\,dt \biggr|^2 .
\]

Now, as a priori the variance $Q_0(u) \leq
E_{\nu_\rho}[J^2_{-1,0}(u)]
= O(\sqrt{u})$, we need only find the limit of
%
\begin{eqnarray}\label{Q(T)} Q_1(T) & = &
\frac{\rho(1-\rho)}{\sqrt{T}}
\sum_i \biggl| \frac{1}{2}\int_{\epsilon T}^Tp(t,i+1)-p(t,i)\,dt
\biggr|^2\nonumber
\\[-9pt]
\\[-9pt]
&=&\frac{\rho(1-\rho)}{4\sqrt{T}}\sum_i \int_{[\epsilon T,T]^2}
[p(t,i+1)-p(t,i) ]  [p(s,i+1)-p(s,i) ]\,ds\,dt.
\nonumber
\end{eqnarray}

To estimate the integrand, from Doob's inequality, note
%
\begin{eqnarray}
\label{local1}
p(t,x) & =& \mathbb E[P(S_{N_t} = x)] \nonumber\\[-8pt]\\[-8pt]
& =& \mathbb E\Bigl[ P(S_{N_t}=x),
\sup_{t\in[\epsilon T,T]}|N_t/t -1|\leq\epsilon\Bigr] +
O(T^{-10}),\nonumber
\end{eqnarray}
where $N_t$ is a Poisson process with rate $1$ independent of the
discrete time random walk $\{S_k\}$, $N_t/t -1$ is a martingale and
$\mathbb E$
refers to expectation with respect to
$N_t$.
Further (since we could not find an appropriate continuous time
version), from the local limit theorem (Petrov~\cite{Petrov}, Theorem
VII.13; page 205),
uniformly over $x$, with respect to the
discrete time walk, we have for $N_t\geq1$ that
%
\begin{eqnarray}
\label{local2}
P(S_{N_t}=x) & = &
\frac{1}{\sqrt{2\pi N_t}}e^{-\fracc{x^2}{2N_t}}+
\frac{1}{\sqrt{2\pi}}e^{-\fracc{x^2}{2N_t}}\frac{q_2(x/\sqrt{N_t})}{
N_t^{3/2}} +o(N^{-3/2}_t)\nonumber\\[-8pt]\\[-8pt]
&=&\frac{1}{\sqrt{2\pi N_t}}e^{-\fracc{x^2}{2N_t}} +
O(N^{-3/2}_t),\nonumber
\end{eqnarray}
where $q_2(y)=(\gamma_4/24\theta^4)(y^4-6y^2+3)$, $\gamma_k$ is the
$k$th order cumulant and $\theta^2$ is the variance of the symmetric
Bernoulli variable. [In our case, in Petrov's formula,
$q_1(y) = (\gamma_3/6\theta^3)(y^3-3y)\equiv0$ as $\gamma_3=0$.]

Let $p^N(t,x) = P(S_{N_t}=x)$ and $p^R(s,x) =P(S_{R_s}=x)$
where $R_s$ is an independent Poisson process also with rate $1$.
We now argue that only the leading terms in (\ref{local1}) and (\ref{local2})
are significant.

Since $\sum_x p(u,x), \sum_x p^N(u,x)\leq1$,
the error term on order $O(T^{-10})$ in
(\ref{local1}) can be neglected in estimating (\ref{Q(T)}). Indeed,
\[
\frac{O(T^{-10})}{\sqrt{T}}\int_{[\epsilon T, T]^2}\sum_ip(s,i)\,ds\,dt
= \frac{O(T^{-10})}{\sqrt{T}}\int_{[\epsilon T, T]^2}\sum_ip^N(t,i)\,ds\,dt
= o(1).
\]

Also, note the error term of order $O(N^{-3/2}_t)$ in (\ref{local2})
is not
significant with respect to (\ref{Q(T)}). Indeed,
\begin{eqnarray*}
&&\sum_x \frac{1}{\sqrt{2\pi N_t}}\bigl|e^{-\fracc{(x+1)^2}{2N_t}} -
e^{-\fracc{x^2}{2N_t}}\bigr| \\
&& \qquad = \sum_x \frac{1}{\sqrt{2\pi
N_t}}\bigl|e^{-(2x+1)/2N_t}-1\bigr|e^{-\fracc{x^2}{2N_t}}\\
&&\qquad\leq Ce^{-\sqrt{N_t}/4} + C\sum_{|x|\leq N_t^{3/4}}\frac
{1}{\sqrt{2\pi
N_t}}\frac{|x|}{N_t} e^{-x^2/2N_t}\\
&&\qquad \leq\frac{C}{\sqrt{N_t}}
\end{eqnarray*}
for some constants $C$. Then,
given $|N_t/t-1|,|R_s/s-1|\leq\epsilon$ for $s,t\in
[\epsilon T,T]$,
a product of $\sum_i (2\pi N_t)^{-1/2}|e^{-(i+1)^2/2N_t} -
e^{-i^2/2N_t}|$
and the error term with respect to the $s$-integration, for instance,
leads to bounding
\begin{eqnarray*}
&&\frac{1}{\sqrt{T}}\int_{[\epsilon T,T]^2}
\sum_i\frac{1}{R_s^{3/2}}\frac{1}{\sqrt{2\pi N_t}}\bigl|e^{-\fracc
{(i+1)^2}{2N_t}} -
e^{-\fracc{i^2}{2N_t}}\bigr|\,ds\,dt\\
&&\qquad\leq\frac{1}{\sqrt{T}}\int_{[\epsilon
T,T]^2}\frac{C}{R_s^{3/2}\sqrt{N_t}}\,ds\,dt\leq O(T^{-1/2}).
\end{eqnarray*}

Therefore, focusing on the leading order terms,
\begin{eqnarray*}
&&\frac{1}{4\sqrt{T}}\sum_i \int_{[\epsilon T,T]^2}
[p^N(t,i+1) - p^N(t,i) ] [p^R(s,i+1)-p^R(s,i) ]\,ds\,dt \\
&&\qquad= o(1) +
\frac{1}{8\pi\sqrt{T}}\sum_i \int_{[\epsilon T,T]^2}
\frac{1}{\sqrt{R_sN_t}} \bigl[e^{-(i+1)^2/2N_t} - e^{-i^2/2N_t} \bigr]\\
&&\hphantom{o(1) +
\frac{1}{8\pi\sqrt{T}}\sum_i \int_{[\epsilon T,T]^2}}\qquad\quad{}\times
\bigl[e^{-(i+1)^2/2R_s}-e^{-i^2/2R_s} \bigr]\,ds\,dt.
\end{eqnarray*}
Now, using again $|N_t/t-1|, |R_s/s-1|\leq\epsilon$ for
$s,t\in[\epsilon T,T]$, we further evaluate the integral on the
right-hand side as
\begin{eqnarray*}
&&o(1)+\frac{1}{8\pi\sqrt{T}}\sum_{|i|\leq T^{3/4}} \int
_{[\epsilon T,T]^2} \frac{i^2}{R_sN_t\sqrt{R_sN_t}}e^{\fracd{-i^2}{2}
[\fraca{1}{N_t}+\fraca{1}{R_s} ]}\,ds\,dt\\
&&\qquad= o(1)+\frac{1}{8\pi\sqrt{T}}\int_{[\epsilon T,T]^2}\int
_{-\infty}^{\infty} \frac{x^2}{R_sN_t\sqrt{R_sN_t}}e^{\fracd{-x^2}{2} [\fraca{1}{N_t}+\fraca{1}{R_s} ]}\,dx\,ds\,dt\\
&&\qquad= o(1)+\frac{\sqrt{2}}{8\sqrt{\pi T}}\int_{[\epsilon T,T]^2}
(N_t+R_s)^{-3/2}\,ds\,dt =: Q_2(T,\epsilon).
\end{eqnarray*}
Finally,
we have that $Q_2(T,\epsilon)$ satisfies
\[
\lim_{T\uparrow\infty}\biggl| Q_2(T,\epsilon) - \frac{\sqrt
{2}-1}{\sqrt{\pi}}\biggr| \leq c(\epsilon),
\]
where $c(\epsilon)$ vanishes
as $\epsilon\downarrow0$.
\end{pf}

\subsection{\texorpdfstring{Proof of Theorem \protect\ref{mainthm4}}{Proof of Theorem 1.8}}
We concentrate on the argument for the tagged particle, as a similar
proof holds for the current.
By symmetry,
\[
P\bigl(|X(N^2T)|/N\geq a\bigr)= 2P\bigl(X(N^2T)/N\geq a\bigr).
\]
From (\ref{positive_a}), and noting
$J_{-1,0}(t)-J_{\lfloor aN\rfloor,\lfloor aN\rfloor+1}(t) =
\sum_{x=0}^{\lfloor aN\rfloor}\eta_t(x)-\eta_0(x)$ by the
development of  Section
\ref{tagged_particle_current_subsection}, we have
\[
\{X(N^2t)\geq aN\} = \Biggl\{J_{\lfloor
aN\rfloor,\lfloor aN\rfloor+1}(N^2 t) \geq\sum_{x=0}^{\lfloor
aN\rfloor}\eta_0(x)\Biggr\}.
\]

We now rewrite currents in terms of the standard Harris stirring
process $\{\xi^x_t\}$. Namely, at time $t=0$, a particle is put at
each $x\in\Z$.
Then, to bonds $(x,x+1)$ in $\Z$, associate independent Poisson clocks
with parameter $1/2$. When the clock rings at a bond, interchange the
positions of the particles at the bond's vertices. Let $\xi^x_t$ be
the position at time $t$ of the particle initially at~$x$. Then the
exclusion process, starting from initial configuration~$\eta$,
satisfies $\eta_t(x) = 1 \{x\in\{\xi^i_t\dvtx \eta(i)=1\} \}$. More
details and constructions can be found in Chapter VIII~\cite{Liggett}.

Then, for $0\leq a\leq1$,
\[
J_{\lfloor
aN\rfloor,\lfloor aN\rfloor+1}(N^2t) = \sum_{x\leq\lfloor
aN\rfloor}
\eta_{0}(x)1_{[\xi^x_{N^2t}>\lfloor aN\rfloor]}
- \sum_{x> \lfloor aN\rfloor}\eta_0(x)1_{[\xi^x_{N^2t}\leq\lfloor
aN\rfloor]}.
\]
Write, given the initial profile $\eta_0$ is deterministic,
by Chebyshev, that
%
\begin{eqnarray}
\label{thm4_1}
\frac{1}{N}\log P\bigl(X(N^2t)\geq aN\bigr) &\leq&
\frac{1}{N}\log E\exp\Biggl\{-\lambda\sum_{x=0}^{\lfloor aN\rfloor
}\eta_0(x)\Biggr\}\nonumber\\[-8pt]\\[-8pt]
&&{} +\frac{1}{N}\log E\exp\bigl\{\lambda J_{\lfloor
aN\rfloor,\lfloor aN\rfloor+1}(N^2 t) \bigr\}.\nonumber
\end{eqnarray}
The first term on the right-hand side tends to $-\lambda a$ as
$N\uparrow\infty$.
The second term is bounded, by Chebyshev and Liggett
\cite{Liggett}, Proposition VIII.1.7, noting $e^{\alpha\sum_{i=k}^l
1_{[x_i\in A]}}$ is positive definite for any $\alpha\in\R$, and
$\log(1+x)\leq x$ for $x\geq1$, by
\begin{eqnarray*}
&&\frac{1}{2N}\log E \exp\biggl\{2\lambda \sum_{x\leq\lfloor
aN\rfloor}
\eta_{0}(x)1_{[\xi^x_{N^2t}>\lfloor aN\rfloor]}\biggr\}\\
&&\quad{} + \frac{1}{2N}\log E \exp\biggl\{-2\lambda\sum_{x>
\lfloor
aN\rfloor}
\eta_{0}(x)1_{[\xi^x_{N^2t}\leq\lfloor aN\rfloor]}\biggr\}\\
&& \qquad\leq \frac{1}{2N}\sum_{x\leq\lfloor
aN\rfloor} \bigl(e^{2\lambda\eta_0(x)}-1\bigr) P( \xi^x_{N^2t}>\lfloor
aN\rfloor)\\
&& \qquad  \quad{} + \frac{1}{2N}\sum_{x> \lfloor
aN\rfloor}\bigl (e^{-2\lambda\eta_0(x)}-1\bigr)P( \xi^x_{N^2t}\leq
\lfloor
aN\rfloor).
\end{eqnarray*}
Given $\eta_0(x) =1_{[|x|\leq N]}$ and $\xi^x_{N^2t}$ marginally is
the position of a simple random walk, started at $x$ at time $N^2t$;
as $N\uparrow\infty$, we have
\[
 \frac{1}{2N}\sum_{x\leq\lfloor
aN\rfloor} \bigl(e^{2\lambda\eta_0(x)}-1\bigr)P( \xi^x_{N^2t}>\lfloor
aN\rfloor)  \rightarrow \frac{e^{2\lambda}
-1}{2}\int_{-1}^{a}P\bigl(N(0,t)> a-x\bigr)\,dx
\]
and
\begin{eqnarray*}&&\frac{1}{2N}\sum_{x> \lfloor
aN\rfloor} \bigl(e^{-2\lambda\eta_0(x)}-1\bigr)P( \xi^x_{N^2t}\leq
\lfloor
aN\rfloor) \\
&&\qquad\rightarrow \frac{e^{-2\lambda}-1}{2}\int_a^1
P\bigl(N(0,t)\leq a -x\bigr)\,dx,
\end{eqnarray*}
where $N(0,t)$ is a normal distribution with mean $0$ and variance
$t$.

Hence, combining the estimates, we have that (\ref{thm4_1}) is less than
\[
-\lambda a +\frac{e^{2\lambda}
-1}{2}\int_{-1}^{a}P\bigl(N(0,t)> a-x\bigr)\,dx
+ \frac{e^{-2\lambda}-1}{2}\int_a^1
P\bigl(N(0,t)\leq a -x\bigr)\,dx.
\]
Choosing $\lambda= \epsilon a$ for small $\epsilon>0$, we obtain
further that (\ref{thm4_1}) is bounded by
\begin{eqnarray*}
&& -\epsilon a^2\biggl[ 1 - \frac{1}{a}\int_{1-a}^{1+a}
P\bigl(N(0,t)>y\bigr)\,dy\biggr] +
O(\epsilon^2a^2) \leq- Ca^2
\end{eqnarray*}
for a constant $C$, noting $1>a^{-1}\int_{1-a}^{1+a}P(N(0,t)>y)\,dy$ for $0<
a\leq1$.

For $a\geq1$, we write
\[
J_{\lfloor aN\rfloor,\lfloor aN\rfloor+1}(t) = \sum_{|x|\leq N} \eta
_0(x)1_{[\xi^x_{N^2t}> \lfloor aN\rfloor]}.
\]
Then, as above,
\begin{eqnarray*}
P\bigl(X(N^2t)\geq aN\bigr) &\leq& e^{-\lambda\sum_{x=0}^{\lfloor
aN\rfloor}\eta_0(x)}E \exp\biggl\{\lambda\sum_{|x|\leq N}\eta
_0(x)1_{[\xi^x_{N^2t}> \lfloor aN\rfloor]}\biggr\}\\
&\leq& e^{-\lambda N}\prod_{|x|\leq N} E\exp\bigl\{\lambda1_{[\xi
^x_{N^2t}> \lfloor aN\rfloor]} \bigr\}.
\end{eqnarray*}
Taking the logarithm, dividing by $N$ and taking the limit, we obtain
\begin{eqnarray*}
&&\limsup_{N\uparrow\infty} \frac{1}{N}\log P\bigl(X(N^2t)\geq
aN\bigr) \\
&&\qquad\leq-\lambda+ \limsup_{N\uparrow\infty}
\frac{1}{N}\sum_{|x|\leq N} (e^{\lambda}-1) P(\xi^x_{N^2t}> \lfloor
aN\rfloor)\\
&&\qquad\leq -\lambda+ (e^{\lambda}-1)\int_{-1}^1 P\bigl(N(0,t)\geq a-x\bigr)\,dx.
\end{eqnarray*}

Optimizing on $\lambda$, the right-hand side of the above display is
bounded by
\begin{eqnarray*}
\log\int_{-1}^1P\bigl(N(0,t)>a-x\bigr)\,dx + 1-\int_{-1}^1 P\bigl(N(0,t)>a-x\bigr)\,dx < 0.
\end{eqnarray*}
However, for $a$ large, this expression is bounded by $-Ca^2$.

Working with the $0\leq a\leq1$ and $a>1$ bounds, we obtain the desired
quadratic order estimate.

\section{Proofs of approximations}
\label{appendix}
We give the proofs of Propositions~\ref{LY_prop} and~\ref{claudio},
and Lemmas
\ref{cor_limit} and~\ref{unif_limit}.

\subsection{\texorpdfstring{Proofs of Propositions \protect\ref{LY_prop} and \protect\ref{claudio}}
{Proofs of Propositions 1.3 and 2.1}}
The proofs are through
a series of lemmas
inspired by the scheme in~\cite{Lan} (see also Oelschl\"ager \cite
{Oel}, and Bertini, Landim and Mourragui~\cite{BLM}). As several of the
steps are different, we give some details.

To this end, let $\mu$ be a density such that
$I_0(\mu)<\infty$.
The first lemma states that finite rate densities $\mu$, when
integrated against smooth test functions, are
uniformly continuous in time; cf. Lemma 4.4~\cite{BLM}.
%
\begin{lemma}\label{app_cont}
Let $\eta\in D([0,T]; M_1)$ be a density such that $I_0(\eta)<\infty
$, and let
$\mathfrak J\in C^2_K(\R)$. Then, $s\mapsto
\langle\eta_s, \mathfrak J\rangle= \int\mathfrak J(x)\eta_s(x)\,dx$
is a uniformly
continuous function.
\end{lemma}

\begin{pf}
Let $G\in C^{1,2}_K([0,T]\times\R)$.
As $I_0(\eta)<\infty$, from (\ref{weak_eq}), we infer
\[
l^2(\eta;G) \leq2I_0(\eta)\int_0^T\int
G_x^2(t,x)\eta_t(x)\bigl(1-\eta_t(x)\bigr)\,dx\,dt.
\]
Let $F^\delta$ be a smooth approximations of the indicator $1_{[s,t]}(u)$.
Then, by applying the previous inequality with $G^\delta= F^\delta
\mathfrak J$, we obtain
\begin{eqnarray*}
 \biggl| \int\eta_t \mathfrak J \,dx - \int\eta_s \mathfrak J \,dx
\biggr| &=& \lim_{\delta\downarrow0}\biggl\{ l(\eta;G^\delta) +
\frac{1}{2}\int_0^T\int G^\delta_{xx}\eta_u(x) \,dx\,du\biggr\} \\
  &\leq&|t-s|\|\mathfrak J''\|_{L^1}
+ \sqrt{2I_0(\eta)}|t-s|^{1/2}\|\mathfrak J'\|_{L^2},
\end{eqnarray*}
completing the proof.
\end{pf}

For the remainder of the subsection, let $\mu\in D([0,T];M_1)$ be a
density with finite rate, $I_0(\mu)<\infty$. We now build a
succession of approximating densities in the next lemmas with special
properties.

\begin{lemma}\label{app0}
For each $\epsilon>0$,
there exists a density
$\hat{\mu}$, smooth in the space variable, such that: (1) the
Skorohod distance $d(\hat\mu;\mu)<\epsilon$; (2) there is
$0<\delta_\epsilon<1$ such that $\delta_\epsilon<
\hat{\mu}(t,x)
< 1-\delta_\epsilon$ for $(t,x)\in[0,T]\times\R$; (3)
$|I_0(\hat{\mu})-I_0(\mu)|<\epsilon$.

In addition, (4) if $\hat\gamma\in M_1(\rho_*,\rho^*)$ is piecewise
continuous, $0<\hat\gamma(x)<1$ for $x\in\R$, and $h(\mu_0;\hat
\gamma)<\infty$, then also
$|h(\hat{\mu}_0;\hat\gamma)-h(\mu_0;\hat\gamma)|<\epsilon$.
\end{lemma}

\begin{pf}
For $0<\rho_*,\rho^*<1$, let
$\gamma\in M_1(\rho_*,\rho^*)$
be a function. Consider
%
\begin{equation}
\label{initial_mu}
\mu^{b,\alpha} = \sigma_{t+\alpha}*\gamma+
b(\sigma_\alpha*\mu-\sigma_{t+\alpha}*\gamma)
\end{equation}
for $0\leq b\leq1$ and
$\alpha\geq0$. Clearly, $\mu^{b,\alpha}$ is smooth in the space
variable when $\alpha>0$.\vadjust{\goodbreak}

Next, for fixed $\alpha>0$ and $0<b<1$, there is $0<\delta<1$ such that
$\delta<\mu^{b,\alpha}<1-\delta$ as
$\sigma_{t+\alpha}*\gamma$ is strictly bounded between $0$ and $1$ for
$t\in[0,T]$.

Now,
$\mu^{b,\alpha} \rightarrow\mu^{1,\alpha}$ as $b\uparrow1$ in
$D([0,T]\times
\R)$. Also, noting $\lim_{\alpha\downarrow0}
\|\sigma_\alpha*G - G\|_{L^1(\R)}=0$ for $G\in L^1(\R)$, we have
also have the Skorohod convergence
$\mu^{1,\alpha} \rightarrow\mu$ as $\alpha\downarrow0$.

By lower semi-continuity of $I_0$,
\[
\liminf I_0(\mu^{b,\alpha}) \geq I_0(\mu).
\]
On the other hand, by convexity of $I_0(\nu)$, we have
\begin{eqnarray*}
I_0(\mu^{b,\alpha}) & \leq&
(1-b)I_0(\sigma_{t+\alpha}*\gamma) +
bI_0(\sigma_\alpha*\mu).
\end{eqnarray*}
Note that $I_0(\sigma_{t+\alpha}*\gamma) = 0$, and by
translation-invariance
and convexity,
the right-hand side in the display is less than
\[
b\int\sigma_\alpha(y)I_0\bigl(\mu(t,x-y)\bigr)\,dy = bI_0(\mu)\uparrow1
\qquad\mbox{as } b\uparrow1.
\]

Similarly, if $\hat\gamma\in M_1(\rho_*,\rho^*)$ is piecewise
continuous, $0<\hat\gamma<1$ and $h(\mu_0;\hat\gamma)<\infty$,
then, by lower semi-continuity and
convexity of
$h(\cdot;\hat\gamma),$ we have
\[
h(\mu_0;\hat\gamma)\leq\liminf_{b\uparrow1, \alpha\downarrow0}
h(\mu^{b,\alpha}_0;\hat\gamma)
\]
and
\[
h(\mu^{b,\alpha}_0;\hat\gamma)\leq (1-b) h(\sigma_\alpha*\gamma
;\hat\gamma) +
b h(\sigma_\alpha*\mu_0;\hat\gamma).
\]
Also, once more by convexity,
\[
h(\sigma_\alpha*\mu_0;\hat\gamma) \leq \int dy\sigma_\alpha
(y)\int dx
h_d\bigl(\mu_0(x);\hat\gamma(x-y)\bigr).
\]
The right-hand side, since $|h(\mu_0(\cdot); \hat\gamma(\cdot- y))
- h(\mu_0;\hat\gamma)| \leq
C|y|$
by properties of $\hat\gamma$,
converges to $h(\mu_0;\hat\gamma)$ as $\alpha\downarrow0$.
By
the same argument, $\lim_{\alpha\downarrow0}h(\sigma_\alpha*\gamma
;\hat\gamma)=h(\gamma;\hat\gamma)$. Hence $\lim_{b\uparrow1,
\alpha\downarrow0} h(\mu^{b,\alpha}_0;\hat\gamma) = h(\mu_0;\hat
\gamma)$.

Therefore, statements (1)--(4) hold for $\hat{\mu} = \mu^{b,\alpha
}$ when $b\sim1$, $\alpha\sim0$.
\end{pf}

\begin{lemma}\label{app1} Let $\hat\mu$ be the density
constructed in Lemma
\ref{app0}.
Then: (1) for each $\epsilon>0$,
there exists a smooth density
$\tilde\mu$ such that $\tilde\mu_0=\hat\mu_0$; (2) the Skorohod
distance $d(\tilde\mu;\hat\mu)<\epsilon$; (3)
$|I_0(\tilde{\mu})-I_0(\hat\mu)|<\epsilon$.
Also, (4) all partial derivatives of $\tilde\mu$ are uniformly
bounded in $[0,T]\times\R$.
\end{lemma}

\begin{pf} To obtain a smooth density, we need only approximate
$\hat\mu$ by smoothing in the time variable.
Define for
$\beta>0$ a density which is constant in time on a short time interval.
\[
\nu^{\beta}(t,x) = \cases{
\hat\mu_0(x), &\quad for $0\leq t<
\beta$,\cr
\hat\mu(t-\beta,x),&\quad for $\beta\leq t\leq T+\beta$.
}
\]
Let $\kappa_\varepsilon\in C^\infty_K(\R)$ be smooth approximations
of the identity in $L^1(\R)$ such that\vadjust{\goodbreak} $\kappa_\varepsilon\geq0$,
$\int
\kappa_\varepsilon(x)\,dx = 1$, $\operatorname{Supp}(\kappa
_\varepsilon) \subset
(0,\varepsilon)$ and for $f\in L^1(\R)$, $f*\kappa_\varepsilon
\rightarrow f$ as $\varepsilon\downarrow0$ in $L^1$.
Form the convolution, for $0<\varepsilon\leq\beta$,
\[
\nu^{\beta,\varepsilon}(t,x) = \int_0^T \nu^{\beta}(t+s,x)\kappa
_\varepsilon(s)\,ds.
\]
It is clear, by continuity of $\hat\mu$ in time (Lemma~\ref{app_cont}),
that
$\lim_{\beta\downarrow0}\lim_{\varepsilon\downarrow
0}\nu^{\beta,\varepsilon} = \hat\mu$ in
$D([0,T]; M_1)$.
By construction, $\nu^{\beta,\varepsilon}$ is smooth, and
also $\nu^{\beta,\varepsilon}_0 = \hat\mu_0$.

From lower semi-continuity and convexity
\[
\liminf_{\beta,\varepsilon}
I_0(\nu^{\beta,\varepsilon}) \geq
I_0(\hat\mu) \quad\mbox{and}\quad
I_0(\nu^{\beta,\varepsilon})
\leq\int_0^T
\kappa_\varepsilon(s)I_0\bigl(\nu^{\beta}(t+s,x)\bigr)\,ds.
\]
Using the variational definition of $I_0$,
noting $\hat\mu_0 = \sigma_\alpha*(\gamma+ b(\mu_0-\gamma))$,
the rate of $\nu^\beta$ on the interval $[0,\beta]$ is bounded by
\begin{eqnarray*}
\sup_{G\in C^{1,2}_K} \frac{1}{2}\int_0^\beta\int G_x \partial
_x\hat\mu_0 - G_x^2\hat\mu_0(1-\hat\mu_0)\,dx\,dt
& \leq&
\frac{\beta}{8}\int
\frac{(\partial_x\hat\mu_0)^2}{\hat\mu_0(1-\hat\mu_0)}\,dx,
\end{eqnarray*}
which vanishes as $\beta\downarrow0$.
On the other hand, by formula (\ref{I_0}),
the rate of $\nu^\beta$ on the interval
$[\beta,T]$ converges to $I_0(\hat\mu)$ as $\beta\downarrow0$. We
can conclude then that $\lim_{\beta,\varepsilon\downarrow0}
I_0(\nu^{\beta,\varepsilon}) \rightarrow I_0(\hat\mu)$.

Moreover, by differentiating the convolutions, since $\|\nu^{\beta
,\varepsilon}\|_{L^\infty}\leq1$,
we have $\|\partial^{(k)}_x\partial^{(l)}_t\tilde
\mu\|_{L^\infty}\leq\|\partial_x^{(k)}\sigma_\alpha\|_{L^1}\|
\partial_t^{(l)}\kappa_\varepsilon\|_{L^1}<\infty$.

Hence, to find the desired density, we can take $\tilde\mu= \nu
^{\beta,\varepsilon}$ for $\beta,\varepsilon$ small.
\end{pf}

We now continue to adjust the approximation so that the associated
function ``$H_x$'' of the approximating density has desired properties.

\begin{lemma}\label{app2} Let $\tilde{\mu}$ be the density constructed
in Lemma~\ref{app1}, and $\tilde H_x$ be associated to it via (\ref
{weak_asym_eq}).
Then: (1) $\tilde H_x\in C^\infty([0,T]\times\R)$; (2) $\|\tilde
H_x\|_{L^2([0,T]\times
\R)}<\infty$; (3) $\|\tilde H_x\|_{L^\infty([0,T]\times\R)}<\infty$.
\end{lemma}

\begin{pf} By construction, we recall, for a $\delta>0$,
that $\delta<\tilde\mu<1-\delta$, $\tilde\mu$
is smooth with uniformly bounded derivatives
on $[0,T]\times\R$ of all orders, and
%
\begin{equation}
\label{app_eq_1}\partial_t \tilde{\mu} = \tfrac{1}{2}\partial
_{xx}\tilde{\mu} -
\partial_x [\tilde H_x\tilde\mu(1-\tilde\mu) ].
\end{equation}
Then, as $2I_0(\tilde\mu) = \int_0^T\int
(\tilde H_x)^2\tilde\mu(1-\tilde\mu)\,dx\,dt<\infty$, we obtain the
$L^2$ bound
on $\tilde H_x$, and, by solving for $\tilde H_x$ in
(\ref{app_eq_1}), we have that $\tilde H_x$ is smooth.

We now deduce that $\tilde H_x$ is bounded in $L^\infty$. This bound
will follow from the $L^2$ bound on $\tilde H_x\tilde\mu(1-\tilde\mu
)$ and $\delta<\tilde\mu<1-\delta$, if we show that $\tilde
H_x\tilde\mu(1-\tilde\mu)$ is Lipschitz in both space and time
variables with uniform constant over $[0,T]\times\R$. However, from
(\ref{app_eq_1}),
$\tilde H_x\hat\mu(1-\hat\mu)$ is Lipschtiz in the space variable
with uniform constant as $\partial_t\tilde
\mu$
and $\partial_{xx}\tilde\mu$ are bounded on $[0,T]\times\R$.

To show $\tilde H_x\tilde\mu(1-\tilde\mu)$ is also Lipschitz in the
time variable $t$ uniformly over $[0,T]\times\R$, write
\begin{eqnarray*}
&&\tilde H_x\tilde\mu(1-\tilde\mu)(x,t) \\
&&\qquad= \tilde H_x\tilde\mu(1-\tilde\mu)(0,t)
+ (1/2)\partial_x \tilde\mu(x,t)-(1/2)\partial_x\tilde\mu(0,t) -
\int_0^x\partial_t\tilde\mu \,dy.
\end{eqnarray*}
The first three terms on the right-hand side are clearly uniformly
Lipschitz in $t$ as their partial derivatives in time are bounded on $[0,T]$.

To treat the last term, consider a smooth $G$ compactly supported in
$[-\epsilon,x+\epsilon]$ which equals $1$ on $[0,x]$. Since $\partial
_{tt}\tilde\mu$ is bounded, we have
\begin{eqnarray*}
\biggl|\int_0^x\partial_{tt}\tilde\mu(u,y)\,dy \biggr| &\leq& \biggl|\int G(y)
\partial_{tt}\tilde\mu(u,y)\,dy \biggr| + 2C\epsilon.
\end{eqnarray*}
Now, by construction in the proof of Lemma~\ref{app1}, $\tilde\mu=
\kappa_\varepsilon*\nu^\beta$, and so
\[
\int G(y)\partial_{tt}\tilde\mu(u,y)\,dy = \int_0^T \int
G(y)\kappa''_\varepsilon(s)\nu^{\beta}(u+s,y)\,dy\,ds.
\]
As $I_0(\nu^{\beta})<\infty$, we can associate via
(\ref{weak_asym_eq}) an $H_x^{\beta}$ to the density $\nu^{\beta}$.
From the weak formulation (\ref{weak_eq}), and $\kappa'_\varepsilon
(0)=\kappa'_\varepsilon(T)=0$, the right-hand side equals
\begin{eqnarray*}
&&-\frac{1}{2}\int_0^T \int G''(y)\kappa'_\varepsilon(s)\nu^{\beta
}(u+s,y)\,dy\,ds,\\
&&- \int_0^T\int G'(y)\kappa'_\varepsilon(s)H_x^{\beta}\nu^{\beta
}(1-\nu^{\beta})(u+s,y)\,dy\,ds.
\end{eqnarray*}
The first integral, because $\nu^{\beta}$ is bounded and $G'\neq0$
on a set of width at most $2\epsilon$ is uniformly bounded in time $u$
and space $x$. Similarly, the second integral,
as $\|H_x^{\beta}\nu^{\beta}(1-\nu^{\beta})\|_{L^2}\leq2I_0(\nu
^\beta)<\infty$, is also both uniformly bounded
in $u$ and $x$.
\end{pf}

The function $\tilde H_x$ associated to $\tilde\mu$ in Lemma \ref
{app2}, although smooth,
does not necessarily have compact support. Let
$H_x^m\in C^\infty_K((0,T]\times\R)$ be smooth approximations of
$\tilde H_x$ with the
following properties:
\[
\|H_x^m -\tilde H_x\|_{L^2([0,T]\times\R)} \leq
m^{-1} \quad \mbox{and}\quad\sup_{\stackrel{t\in[0,T]}{x\in
[-m,m]}} | H^m_x - \tilde H_x| \leq m^{-1}.
\]

Denote $w^m\in D([0,T]; M_1)$ as the smooth density with initial
condition
$w^m_0=\tilde\mu_0$ which satisfies the equation
\[
\partial_t v =(1/2)\partial_{xx}v - \partial_x\bigl(H^m_x v(1-v)\bigr).
\]
Existence, for instance, follows from the hydrodynamic limit for
weakly asymmetric exclusion processes in~\cite{KOV} using the
replacement estimates Theorem 6.1\vadjust{\goodbreak} and Claims 1, 2~\cite{LY}, Section 6;
see also Theorem 3.1~\cite{Lan}. Uniqueness in the class
of bounded solutions follows by the method of Proposition~3.5~\cite{Oel}.

We now show that $w^m$, whose associated function $H^m_x\in C_K^\infty
$, is close to $\tilde\mu$. Hence, $w^m$ will turn out to be a
suitable candidate with respect to Proposition~\ref{claudio}. In
addition, we will be able to deduce that $\mu\in\mathcal A$ under
(LEM) initial distributions.

\begin{lemma}\label{app3} The sequence $w^m$
converges uniformly to $\tilde\mu$ on
compact subsets
of $[0,T]\times\R$, and hence in $D([0,T]; M_1)$. Also,
$I_0(w^m)\rightarrow I_0(\tilde\mu)$.
\end{lemma}

\begin{pf} Suppose that we have proven
$w^m
\rightarrow\tilde\mu$ uniformly on compact subsets. As $\|H_x^m-\hat
H_x\|_{L^2}\rightarrow0$, we would then conclude
$I_0(w^m)
\rightarrow I_0(\tilde\mu)$.
In the following, the constant $C$ may change line to line.

Now, given
$\partial_t\sigma_t(x) = (1/2)\partial_{xx}\sigma_t(x)$, we have
for $t,h>0$ that
\begin{eqnarray*}
&&\sigma_h* w^m_t(y) - \sigma_{t+h}*w^m_0(y) \\
&&\qquad= \int_0^t \int H^m_xw^m(1-w^m)(s,z) \frac
{-(z-y)}{t+h-s}\sigma_{t-s+h}(z-y)\,dz\,ds.
\end{eqnarray*}
By properties of $w^m, H_x^m$ and
$(|z-y|/\sqrt{u})\exp(-(z-y)^2/4u) \leq1$,
\begin{eqnarray*}
&&|H_x^mw^m(1-w^m)|(s,z)\frac{|z-y|}{t+h-s}\sigma_{t+h-s}(z-y)\\
&&\qquad\leq
C|t-s|^{-1/2}\sigma_{2T}(z-y)\in L^1([0,t]\times\R).
\end{eqnarray*}
Hence, taking $h\downarrow0$, we obtain
%
\begin{eqnarray}\label{app_eq_2}
w^m_t(y) & = & \sigma_t*w^m_0(y) \nonumber\\[-8pt]\\[-8pt]
&& {}+ \int_0^t \int H^m_xw^m(1-w^m)(s,z)
\frac{-(z-y)}{t-s}\sigma_{t-s}(z-y)\,dz\,ds. \nonumber
\end{eqnarray}
Equation (\ref{app_eq_2}) also holds with respect to $\tilde{\mu}$.

Let now $|y|\leq m/2$. We have then, using again
$(|z-y|/\sqrt{u})\exp(-(z-y)^2/4u) \leq1$, and $w^m_0=\tilde\mu
_0$, that
%
\begin{eqnarray}
\label{app_eq_4}
 &&|w^m_t(y) - \tilde\mu_t(y)|\nonumber\\
  && \qquad \leq\sigma_t*|w^m_0 - \tilde\mu
_0|(y) \nonumber\\
&& \qquad  \quad  {}+ \int_0^t
\int | H^m_xw^m(1-w^m) -
\hat H_x\tilde\mu(1-\tilde\mu) |(s,z)\frac{|z-y|}{t-s}\sigma
_{t-s}(z-y)\,dz\,ds \\
 && \qquad \leq C\int_0^t
\int | H^m_xw^m(1-w^m) -
\tilde H_x\tilde\mu(1-\tilde\mu) |(s,z)\nonumber\\
&& \qquad  \quad \hphantom{C\int_0^t
\int } {}\times(t-s)^{-1/2}\sigma_{2(t-s)}(z-y)\,dz\,ds.
\nonumber
\end{eqnarray}

We now estimate
the last line
in two parts, noting
\begin{eqnarray*}
&& | H^m_xw^m(1-w^m) -
\tilde H_x\tilde\mu(1-\tilde\mu) |(s,z) \\
&&\qquad\leq| H^m_x - \tilde H_x|(s,z) +
|\tilde H_x| |\tilde\mu(1-\hat\mu) - w^m(1-w^m) |(s,z).
\end{eqnarray*}
The first part, noting $|y|\leq m/2$, by properties of $H^m_x$,
$\|\sigma_t(x)\|_{L^2([0,T]\times\R)}\leq CT^{1/4}$, and
$\sup_{t\in
(0,T]}t^{-1/2}\sigma_{4t}(1) \leq C$,
is bounded for large $m$, as
\begin{eqnarray*}
&&\int_0^t\int|H^m_x-\tilde H_x|(s,z) (t-s)^{-1/2}
\sigma_{2(t-s)}(z-y)\,dz\,ds\\
&&\qquad\leq\int_0^t\int_{|z|\geq m} |H^m_x-\tilde H_x|(s,z) (t-s)^{-1/2}
\sigma_{2(t-s)}(z-y)\,dz\,ds + m^{-1}\sqrt{t}\\
&&\qquad\leq C\int_0^t \int_{|z|\geq m}|H^m_x-\tilde H_x|(s,z)\sigma
_{4(t-s)}(z-y)\,dz\,ds +
m^{-1}\sqrt{t}\\
&&\qquad\leq Cm^{-1}T^{1/4}
+ m^{-1}\sqrt{T}.
\end{eqnarray*}
The second part is decomposed as the sum of three terms,
\begin{eqnarray*}
&&\int_0^t\int|\tilde H_x| |\tilde\mu(1-\tilde\mu)-w^m(1-w^m)
|(s,z) (t-s)^{-1/2}
\sigma_{2(t-s)}(z-y)\,dz\,ds \\
&&\qquad= D_1 + D_2 + D_3,
\end{eqnarray*}
where $D_1,D_2,D_3$ is the integral over $[0,t]\times\{|z|\geq
m/2+\epsilon\}$, $[0,t]\times\{m/2\leq|z|\leq m/2+\epsilon\}$ and
$[0,t]\times\{|z|\leq m/2\}$, respectively, for $\epsilon>0$.

The term $D_1$, noting $\sup_{t\in(0,T]}t^{-1/2}\sigma_{4t}(\epsilon
) \leq
C_\epsilon$, is bounded by
\begin{eqnarray*}
&&2\int_0^t\int_{|z|\geq m/2+\epsilon} |\tilde H_x|(s,z) (t-s)^{-1/2}
\sigma_{2(t-s)}(z-y)\,dz\,ds\\
&&\qquad\leq C(\epsilon, T)\|\tilde H_x\|_{L^2([0,T]\times\{z\dvtx
|z|\geq m/2\})}.
\end{eqnarray*}
The second term $D_2$ is bounded by
\begin{eqnarray*}
&&2\int_0^t\int_{m/2\leq|z|\leq m/2+\epsilon} |\tilde
H_x||t-s|^{-1/2}\sigma_{2(t-s)}(z-y)\,dz\,ds\\
&&\qquad\leq C\|\tilde H_x\|_{L^\infty} \int_0^ts^{-3/4}\int
_0^\epsilon
s^{-1/4}e^{-z^2/4s}\,dz\,ds\\
&&\qquad\leq C\|\tilde H_x\|_{L^\infty}T^{1/4}\sqrt{\epsilon}.
\end{eqnarray*}
The third term $D_3$ is bounded, with respect to a $\tau\geq t$, by
%
\begin{eqnarray}
\label{D_3_eq}
&&2\int_0^t\int_{|z|\leq m/2}|\tilde H_x||\tilde\mu-
w^m|(s,z)|t-s|^{-1/2}\sigma_{2(t-s)}(z-y)\,dz\,ds\nonumber\\[-8pt]\\[-8pt]
&&\qquad\leq2\sqrt{t}\|\tilde H_x\|_{L^\infty} \sup_{\stackrel
{|z|\leq m/2}{s\leq
\tau}}|\tilde\mu_s(z)-w^m_s(z)|.\nonumber
\end{eqnarray}
Hence, for $\tau>0$ small enough but fixed, which satisfies $2\sqrt
{\tau}\|\tilde H_x\|_{L^\infty} = 1/2$, or $\tau= (16\|\tilde H_x\|
_{L^\infty})^{-1}$, and $L <m/2$, we have
\begin{eqnarray*}
 \sup_{\stackrel{|z|\leq L}{t\leq\tau}}|\tilde\mu(z,t)-w^m(z,t)|
&\leq&\sup_{\stackrel{|z|\leq m/2}{t\leq\tau}}|\tilde\mu
(z,t)-w^m(z,t)| \\
& \leq& C(T)m^{-1}+ 2C(\epsilon,T)\|\tilde H_x\|
_{L^2([0,T]\times\{z\dvtx
|z|\geq m/2\})}\\
&&{} + 2C\|\tilde
H_x\|_{L^\infty}T^{1/4}\sqrt{\epsilon}.
\end{eqnarray*}
Here, we absorbed the right-hand side of (\ref{D_3_eq}) into the
left-hand side above.

We may repeat the same scheme, starting from time $\tau$, where now
the initial difference (\ref{app_eq_4}) is taken into account:
\begin{eqnarray*}
&&\sup_{|y|\leq m/3} \sigma_t*|w^m_\tau-\tilde\mu_\tau|(y) \\
&&\qquad\leq\sup_{|z|\leq m/2}|w^m_\tau- \tilde\mu_\tau|(z)
+ \sup_{|y|\leq m/3}\int_{|z|>m/2}\sigma_t(y-z)\,dz\\
&&\qquad\leq \sup_{|z|\leq m/2}|w^m_\tau- \tilde\mu_\tau|(z)
+ e^{-Cm^2/T}.
\end{eqnarray*}

With a finite number of iterations of such type, say $N_\tau=[T/\tau
]+1$ iterations, when $L< m/N_\tau$,
we obtain uniform convergence, as $m\uparrow\infty$, for $|z|\leq L$
and $0\leq s\leq T$.
\end{pf}

\begin{pf*}{Proof of Proposition~\ref{claudio}} The proof
follows
by applying Lemmas~\ref{app0}, \ref{app1} and~\ref{app3} to build a
density $\mu^+=w^m$, which satisfies specifications (i)--(viii). We
remark property (v) is shown as follows: When
$\mu_0= \gamma$, by construction in (\ref{initial_mu}), we have
$w^m_0=\tilde\mu_0=\hat\mu_0=\mu^{b,\alpha}_0= \sigma_\alpha
*\gamma$. When $\gamma(x)\equiv\rho$, this reduces to $\tilde\mu
_0(x) \equiv\rho$.
\end{pf*}

Starting under (DIC) initial conditions, however, to prove Proposition
\ref{LY_prop}, we will need to specify that $w^m$ can be approximated
by a suitable density with initial value equal to $\mu_0=\gamma\in
M_1(\rho_*,\rho^*)$.
%
\begin{lemma}
\label{app_initial_cond}
Recall $w^m$ from Lemma~\ref{app3}.
Suppose $\mu_0 = \gamma\in M_1(\rho_*,\rho^*)$. Then, for $\epsilon
>0$, $\exists M$ such
that $\forall m\geq M$, there is a density
\mbox{$\bar\chi\in C^\infty((0,T]\times\R)$}, such that: (1) equation
(\ref{weak_asym_eq}) is satisfied with respect to $\bar H_x\in
C^\infty_K([0,T]\times\R)$; (2) initial value $\bar\chi_0 = \gamma$;
(3) the Skorohod distance $d(\bar\chi, w^m)<\epsilon$; (4)
$|I_0(\bar\chi)-I_0(w^m)|<\epsilon$.
\end{lemma}

\begin{pf} Consider $w^m_0$ from Lemma~\ref{app3}. From the
assumption $\mu_0=\gamma$, we have
$w^m_0=\tilde\mu_0 = \sigma_\alpha* \gamma$ from (\ref
{initial_mu}). Form the density $\bar\chi$ as follows:
\[
\bar\chi= \cases{
\sigma_t *\gamma, &\quad for $0\leq t\leq\alpha$,\cr
w^m_{t-\alpha}, &\quad for $\alpha\leq t\leq T$.
}
\]
Since $H^m_x$ is supported on a compact subset of $(0,T]\times\R$,
$\bar\chi\in C^\infty((0,T]\times\R)$, and
satisfies (\ref{weak_asym_eq}) with respect to $\bar H_x\in C^\infty
_K([0,T]\times\R)$ given by
\[
\bar H_x = \cases{
0,&\quad for $(t,x)\in[0,\alpha]\times\R$,\cr
H^m_x(t-\alpha,x), &\quad for $(t,x)\in[\alpha,T]\times\R$.
}
\]

Now,
$2I_0(\bar\chi) = \int_0^T\int\bar H_x^2\bar\chi(1-\bar\chi)\,dx\,dt
= \int_0^{T-\alpha}\int(H^m_x)^2w^m(1-w^m)\,dx\,dt$.
Then, the difference
\[
2I_0(\bar\chi)-2I_0(w^m) = \int_{T-\alpha}^T\int(H^m_x)^2w^m(1-w^m)\,dx\,dt.
\]
To estimate the right-hand side, recall from Lemma~\ref{app3} that
$\|H^m_x - \tilde H_x\|_{L^2}\leq
m^{-1}$, and $w^m\rightarrow\tilde\mu$ uniformly on compact
subsets.
Then
\begin{eqnarray*}
&&\int_{T-\alpha}^T\int(H^m_x)^2w^m(1-w^m)\,dx\,dt\\
&&\qquad\leq2\|H^m_x-\tilde H_x\|^2_{L^2([0,T]\times\R)} +
2\int_0^T\int_{|x|\geq L} \tilde H_x^2 \,dx\,dt \\
&&\qquad\quad{}+ 4\int_0^T\int_{|x|\leq L} \tilde H_x^2 |w^m -
\tilde\mu|\,dx\,dt + 2\int_{T-\alpha}^T\int\tilde
H_x^2\tilde\mu(1-\tilde\mu)\,dx\,dt\\
&&\qquad= B_1 + B_2 + B_3 + B_4.
\end{eqnarray*}
Choose $L=L(\tilde H_x)$ large so that $B_2\leq\epsilon/4$, and take
$m=m(\tilde H_x, L)$ large enough so that both $B_1,B_3\leq
\epsilon/4$.

The term $B_4/4$ is the rate of $\tilde\mu$ on the time interval
$[T-\alpha,T]$. Since $\tilde\mu$ and $\tilde H_x$ depend on $\alpha
$, we bound $B_4$ in terms of $H_x$ (which does not depend on $\alpha
$) to show that it is small when $\alpha$ is small. By the
construction of $\tilde\mu$ in Lemma~\ref{app1}, convexity of the
the rate, translation-invariance and that the rate of $\sigma
_{t+\alpha} *\gamma$ vanishes,
we estimate
\begin{eqnarray*}
B_4 & \leq& 2b\int
\sigma_\alpha(z)\int_0^T\kappa_\varepsilon(s)\int_{T-\alpha}^T
\int H_x^2\mu(1-\mu)(t+s-\beta,x-z)\,dx\,dt\,ds\,dz \\
&\leq& 2 \int_{T-2\alpha}^T \int
H_x^2\mu(1-\mu)\,dx\,dt,
\end{eqnarray*}
when $\beta\leq\alpha\leq T-\alpha$. Then, as $I_0(\mu)<\infty$,
$B_4\downarrow0$ as
$\alpha\downarrow0$.

Hence, with $\alpha$
small enough,
there is $M$ so that for $m\geq M$,
we have $|I_0(\bar\chi)-I_0(w^m)|<\epsilon$. Also, by Lemma~\ref{app3},
$I_0(w^m)\leq I_0(\mu)+1$, and so by uniform continuity (Lemma \ref
{app_cont}),
the Skorohod distance
$d(\bar\chi;w^m)<\epsilon$.
\end{pf}

\begin{pf*}{Proof of Proposition~\ref{LY_prop}} Let $\gamma$ be a
profile associated to an (LEM) or (DIC) measure, and let $\mu$ be such
that $I_\gamma(\mu)<\infty$. By successively
applying Lemmas~\ref{app0},~\ref{app1},~\ref{app3} and \ref
{app_initial_cond}, we can approximate $\mu$ by an appropriate density
$\mu^+$ to verify $\mu\in\mathcal A$. Specifically, under an
(LEM) initial measure, when $I_\gamma(\mu) = I^{LE}_\gamma(\mu)$,
$\mu^+=w^m$ in Lemma~\ref{app3} with
appropriate choice of parameters $b, \alpha, \beta, \varepsilon$ and
$m$. Under a (DIC) initial
configuration, when $I_\gamma(\mu)= I^{DC}_\gamma(\mu)$, $\mu
^+=\bar\chi$ in Lemma~\ref{app_initial_cond} again with suitable
parameters.
\end{pf*}

\subsection{\texorpdfstring{Proof of Lemmas \protect\ref{cor_limit}, \protect\ref{unif_limit}}
{Proof of Lemmas 2.2, 2.3}}

We prove the lemmas in succession.

\begin{pf*}{Proof of Lemma~\ref{cor_limit}}
Note that
\begin{eqnarray*}
 |\mu_t(x) -\hat\gamma(x)| &=&\lim_{h\downarrow0} |\sigma_h*(\mu
_t - \hat\gamma)(x)|\\
& \leq&\lim_{h\downarrow0} |\sigma_{t+h}*(\mu_0 - \hat
\gamma)(x)|\\
   &&{} + \int_0^t\int|H_x|\mu(1-\mu)(s,z)|\partial
_z\sigma_{t-s+h}(z-x)|\,dz\,ds.
\end{eqnarray*}
Since $H_x$ has compact support in $[0,T]\times\R$, the second term
on the right-hand side is bounded by
\[
C_H\int_0^t\int_{|z|\leq M_H} \frac{|x-z|}{t-s}\sigma_{t-s}(z-x)\,dz\,ds
\]
for some constants $C_H, M_H$. Since $(|y|/\sqrt{s})e^{-y^2/4s}\leq
1$, when $|x|\geq M_H$, we can bound it further by
$4C_H\sqrt{T}e^{-(x-M_H)^2/8T}$, which vanishes as $|x|\uparrow\infty$.

The first term, however, is bounded as follows:
\begin{eqnarray*}
&&|\sigma_{t+h}*(\mu_0 - \hat\gamma)(x)|\\
&&\qquad\leq\sup_{|z|\leq l} |\mu_0-\hat\gamma|(z-x) \cdot\int
_{|z|\leq l} \sigma_{t}(z)\,dz + \sqrt{2}e^{-l^2/4T} \int_{|z|\geq
l}\sigma_{2t}(z)\,dz\\
&&\qquad\leq\sup_{|z|\leq l}|\mu_0-\hat\gamma|(z-x) + \sqrt
{2}e^{-l^2/4T}.
\end{eqnarray*}
Now, since $h(\mu_0;\hat\gamma)<\infty$, $\hat\gamma\in M_1(\rho
_*,\rho^*)$, $\|\partial_x\mu_0\|_{L^\infty}<\infty$,
we conclude, for fixed $l$, that $\lim_{|x|\uparrow\infty}\sup
_{|z|\leq l}|\mu_0-\hat\gamma|(z-x)=0$. This completes the proof.
\end{pf*}

\begin{pf*}{Proof of Lemma~\ref{unif_limit}}
Consider Hellinger's inequality $(\sqrt{\alpha}-\sqrt{\beta})^2
\leq
h_d(\alpha;\beta)$. [Let $H(\alpha;\beta) =
(\sqrt{\alpha}-\sqrt{\beta})^2 + (\sqrt{1-\alpha} -
\sqrt{1-\beta})^2$. By Jensen's inequality and\vadjust{\goodbreak} $\log(1-x)\leq-x$ for
$0\leq x < 1$,
$h_d(\alpha;\beta) \geq-2\log[1-(1/2)H(\alpha;\break \beta)] \geq
H(\alpha;\beta)$.] We write then
$(\mu_t(x)-\rho)^2 \leq2(\sqrt{\mu}-\sqrt{\rho})^2 \leq2h_d(\mu
_t(x);\rho)$. Hence,
by Proposition~\ref{energyLemma} (with respect to the density on
$[0,t]\times\R$ with $\hat\gamma\equiv\rho$),
\begin{eqnarray*}\int\bigl(\mu_t(x) - \rho\bigr)^2\,dx & \leq& 2\int h_d(\mu
_t(x);\rho)\,dx \\
&\leq& 4\int_0^t\int(H_x^n)^2\mu(1-\mu)\,dx\,ds \leq8I_0(\mu)
\end{eqnarray*}
uniformly in $0\leq t\leq T$.
\end{pf*}

\section*{Acknowledgments}
We would like to thank Claudio Landim,
K. Ravishankar and Ellen Saada for helpful conversations. Thanks
also to theeree for constructive comments with respect to
extensions at the end of Section~\ref{sec0}, and in other parts.

\printaddresses

\end{document}